\documentclass[12pt]{article}

\usepackage[reqno]{amsmath}
\usepackage{amssymb, amsthm}
\usepackage{young, pst-node}
\usepackage{enumerate} 
\usepackage{url}

\newtheoremstyle{plainsl}%
	{\topsep}
	{\topsep}
	{\slshape} 
	{}
	{\normalfont\bfseries}
	{.}
	{ }
	{}

\swapnumbers

\theoremstyle{plainsl}
\newtheorem{theorem}{Theorem}[section]
\newtheorem{lemma}[theorem]{Lemma}
\newtheorem{corollary}[theorem]{Corollary}
\newtheorem{proposition}[theorem]{Proposition}

\newtheorem{eg}[theorem]{Example}

\renewcommand\proof{\noindent\textsl{Proof. }}
\newcommand\sqr[2]{{\vbox{\hrule height.#2pt
    \hbox{\vrule width.#2pt height#1pt \kern#1pt
        \vrule width.#2pt}\hrule height.#2pt}}}

\renewcommand\qed{%
	\ifmmode\eqno\sqr53
	\else\nolinebreak\ \hfill\sqr53\medbreak\fi}

\numberwithin{equation}{section}

\newcommand{\remove}[1]{}
\newcommand\sgn{\mathrm{sgn}}
\newcommand\aut{\mathrm{Aut}}
\newcommand\fix{\mathrm{fix}}
\newcommand\ind[1]{{\mathrm{ind}(1_{#1})} }
\newcommand\indg[2]{{\mathrm{ind}_{#2}(1_{#1})} }

\newcommand\stab[2]{{\mathrm{stab}_{#2}({#1})}}

\newcommand\la{\lambda}

\title{Multiplicity-free permutation representations of the symmetric group}

\author{
Chris Godsil  and Karen Meagher \footnote{Research supported by NSERC.}\\
\small  Department of Combinatorics and Optimization \\[-0.8ex]
\small University of Waterloo,  Waterloo, Ontario, Canada\\[-0.8ex]
\small \texttt{cgodsil@math.uwaterloo.ca}, \small \texttt{kmeagher@math.uwaterloo.ca}
}

\begin{document}
\maketitle

\begin{abstract}
We determine all the multiplicity-free representations of the
symmetric group. This project is motivated by a combinatorial problem
involving systems of set-partitions with a specific pattern of
intersection.
\end{abstract}


\section{Introduction}

In this paper we aim to determine the multiplicity-free permutation
representations of the symmetric groups. This question of which
subgroups of the symmetric group have a multiplicity-free permutation
representation arose from the combinatorial problem of determining the
largest collection of set-partitions with a specific type of
intersection.  We will start by describing the combinatorial problem
and stating our results and then explain how it is related to the
multiplicity-free representations of the symmetric group. All the
calculations performed in this paper are done using the computational
discrete algebra system GAP~\cite{GAP4}.

\subsection{Set-partitions}
Let $k,\ell$ be integers with $n=\ell k$.  A \textsl{uniform $k$-partition
of an $n$-set} is a set-partition of the $n$-set into $k$ classes each
of cardinality $\ell$.  We say two such partitions $P$ and $Q$ are
\textsl{qualitatively independent} if each cell of $P$ and each cell of $Q$ have 
at least one element in common.  For practical reasons, we want large
sets of pairwise qualitatively independent
partitions~\cite{MR2001c:05037, gargano:92, poljak:83, poljak:89}.

The \textsl{qualitative independence graph} $QI(n,k)$ is the graph
whose vertices are the uniform $k$-partitions, where two partitions
are adjacent if they are qualitatively independent. In this setting,
our problem is to find the cliques of maximum size in $QI(n,k)$.  For
$k=2$ the problem is trivial since the graph $QI(n,2)$ is a
clique. For $k$ a prime power, the size of the maximum clique of
$QI(k^2,k)$ is known to be $k+1$ (\cite{karensthesis}, Section 5.4).
Very few bounds on the clique size of $QI(n,k)$ are known for other
values of $n$ and $k$.

One approach is to make use use of a family of graphs that generalises
the graphs $QI(n,k)$.  If $P$ and $Q$ are two uniform $k$-partitions,
their \textsl{meet table} is the $k \times k$ array with the $i,j$
entry $|P_i \cap Q_j|$. Two meet tables are \textsl{isomorphic} if we
can obtain one from the other by permuting its rows and columns.
Different orderings on the classes in the partitions $P$ and $Q$ could
produce different, but isomorphic meet tables.  Note that $P$ and $Q$
are qualitatively independent if and only if all entries of their meet
table are positive.

For each possible meet table $T$, we define a graph whose vertices are
the uniform partitions, where two partitions are adjacent if their
meet table is isomorphic to $T$. For example, if $n=k^2$ and all
entries in $T$ are equal to one, then the corresponding graph is
$QI(k^2,k)$.  Thus we have a family of graphs, one for each meet
table, whose edges partition the edges of the complete graph on the
uniform partitions.  If the adjacency matrices of these graphs
commute, then standard eigenvalue methods can be used to obtain bounds
on the maximum sizes of cliques.

Mathon and Rosa~\cite{MR86e:05068} show that for $n=9$ and $k=3$ that
the family of adjacency matrices of these graphs is commutative.  This
result is used by Godsil and Newman to prove that the graph $QI(9,3)$
is a core~\cite{Godsil:Newman}. It is conjectured
in~\cite{karensthesis} that the graph $QI(k^2,k)$ is a core for all
values of $k$. A first step towards this is to determine if the
adjacency matrices of these graphs for $n=k^2$ commute for all values
of $k$.

How can we decide if these matrices commute?  The key fact is that the
symmetric group $S_n$ acts as a group of permutations on the set of
uniform partitions, and the matrices we have just defined commute with
each other if and only if this permutation representation is
\textsl{multiplicity free}, that is, each irreducible representation
of $S_n$ occurs in it with multiplicity 0 or 1. (See
\cite{MR0183775}, Section 29.)

\subsection{Results}

We have determined all multiplicity-free permutation representations
of the symmetric group.  Saxl~\cite{MR627512} gives a list of all the
possible subgroups of $S_n$ for $n > 18$ which could be
multiplicity-free. Our work extends Saxl's list not only by including
groups acting on fewer than 18 points but also by determining exactly
which groups in his list are multiplicity free. Mark
Wildon~\cite{Wildon} has also determined which of the groups in Saxl's
list are multiplicity free, his focus is on groups acting on at least 66 points.

\begin{theorem}(Saxl~\cite{MR627512})\label{saxl}
Let $n >18$ and $H \leq S_n$.  If $H$ is multiplicity free, then one
of the following holds:
\begin{enumerate}[(a)]
\item $A_{n-k} \times A_{k} \leq H \leq S_{n-k} \times S_{k}$ for some $k$ with $0 \leq k < n/2$;
\item $n=2k$ and $A_{k} \times A_{k} < H \leq S_{k} \wr S_{2}$;
\item $n=2k$ and $H \leq S_{2} \wr S_{k}$ with $[S_{2} \wr S_{k}:H] \leq 4$;
\item $n=2k+1$ and $H$ fixes a point of $[1,n]$ and is one of the subgroups
in (ii) or (iii) on the rest of $[1,n]$;
\item  $A_{n-k} \times G_{k} \leq H \leq S_{n-k} \times G_{k}$ where $k =
5,6$ or 9 and $G_k$ is $AGL(1,5)$, $PGL(2,5)$ or $P\Gamma L(2,8)$ respectively.
\end{enumerate}
\end{theorem}

Bannai~\cite{MR0357559} determined all the maximal subgroups with
permutation rank less than or equal to five of the finite symmetric
groups.  All of these representations are multiplicity free and,
naturally, occur in our list as well.

There are several infinite families of groups whose permutation
representations is multiplicity free and a long list of groups that do
not belong to any of these families.  A complete list of the groups
can be found in the tables in the
appendix. Table~\ref{table:transitive} is a list of all the transitive
primitive subgroups.  Table~\ref{table:imprimitive} contains all the
transitive imprimitive subgroups and the intransitive groups are given
in Tables~\ref{table:intrans} and \ref{table:intrans2}.

We list the infinite families here. Many of the groups in the list are
subgroups of others in the list, we organized the list so a subgroups
are below the groups that contains them.

\begin{theorem}\label{thm:infinitefam}
The table below is a complete list of the infinite families of
subgroups of the symmetric groups whose permutation representation is
multiplicity free.
\begin{center}
\begin{tabular}{llc}
Subgroup & & Symmetric Group \\
\hline
$S_n$ && $S_n$ \\
$A_n$ && $S_n$ \\

$S_\ell \wr S_2$ && $S_{2\ell}$ \\
$(S_\ell\wr S_2)\cap A_{2\ell}$ &for $k\geq 3$& $S_{2\ell}$ \\
$((S_\ell \times S_\ell) \cap A_{2\ell})\rtimes S_2$ &for $k\geq 3$& $S_{2\ell}$ \\
$A_\ell \wr S_2$ & & $S_{2\ell}$ \\

$S_2\wr S_k$  && $S_{2k}$ \\
$(S_2\wr S_k) \cap A_{2k}$ & $k$ odd & $S_{2k}$ \\

$S_k \times S_{n-k}$ & for $2k \leq n$ & $S_n$ \\
$(S_k \times S_{n-k})\cap A_{n}$ & for $2k\leq n$, $(k,n)\neq (2,4)$& $S_n$\\ 
 $A_k \times S_{n-k}$& for $k \neq 2$, $k \neq n-2$ & $S_n$\\
$A_k \times A_{n-k}$ &for $k\geq 3$ and $2k \leq n-2$ & $S_n$\\

$S_k \times AGL(1,5)$& & $S_{k+5}$\\
$A_k \times AGL(1,5)$ &for $k\geq 4 $ & $S_{k+5}$ \\
$(S_k \times (AGL(1,5)) \cap A_{k+5}$& for $k\geq 5$ & $S_{k+5}$\\

$S_k\times PGL(2,5) $& & $S_{6+k}$ \\
$A_k \times PGL(2,5) $ & for $k\geq 3 $ & $S_{k+6}$ \\
$(S_k \times (PGL(2,5)) \cap A_{k+6}$ & for $k\geq 2$& $S_{k+6}$ \\

$S_k\times P\Gamma L(2,8) $ && $S_{k+9}$ \\
$A_k\times P\Gamma L(2,8) $ & for $k\geq 11$ & $S_{k+9}$  \\

$S_1 \times (S_k \wr S_2)$ && $S_{2k+1}$ \\
$(S_1\times (S_k \wr S_2))\cap A_{2k+1}$ &for $k\geq 3$& $S_{2k+1}$\\
$S_1\times ((S_k^2 \cap A_{2k}) \rtimes S_2)$ &for $k\geq 3$ & $S_{2k+1}$\\

$S_1 \times (S_2 \wr S_k)$&& $S_{2k+1}$ \\
\end{tabular}
\end{center}
\end{theorem}

There are only two infinite family of wreath products of the form
$S_\ell \wr S_k$ whose permutation representation is multiplicity
free: the products $S_\ell \wr S_2$ and $S_2 \wr S_k$.  Except for the
case where $\ell =2$, or $k=2$ and a handful of specific values of
$\ell$ and $k$, the permutation representation of $S_\ell \wr S_k$ is
not multiplicity free. Thus the set of adjacency matrices for the
graphs defined from the meet tables of partitions do not in general
commute.

If $G \leq S_n$ is not one of the groups listed in
Theorem~\ref{thm:infinitefam} and the index of $G$ in $S_n$ is greater
than $\frac{64!}{2^{31}32!}$, then the permutation representation of
$S_n$ on the set of cosets of $G$ is not multiplicity free.

\subsection{Representation Theory}

We will restate basic results from the theory of representations.
Proofs for these results can be found in Chapters 1--4
of~\cite{MR1153249}.

Let $n$ be a positive integer. A partition of $n$ is a set 
$\la = [\la_1,\la_2,\dots ,\la_k]$ where 
\[
\la_i \geq \la_{i+1}, \qquad i \in [1\ldots k-1]
\]
and
\[
\la_1 + \la_2 + \cdots + \la_k =n.
\]
We use $\vdash$ to denote the relation ``is a partition of''; thus $\la\vdash n$
means that $\la$ is a partition of $n$.  If an integer $\la_i$ occurs exactly
$a_i$ times in a partition, we write the partition as
$[\la_1^{a_1},\dots , \la_{\ell}^{a_{\ell}}]$. The conjugate of a
partition $\la$ will be denoted by $\la '$.

It is well-known that for every integer $n$, each of the irreducible
representations of $S_n$ corresponds to a unique partition of $n$.  An
irreducible representation of symmetric group and its corresponding
character will be both denoted by the corresponding partition.

For any group $G$, the trivial representation of $G$ is written as
$1_G$.  For groups $H$ and $G$ with $H \leq G$, if $V$ is a
representation of $H$ then the representation of $G$ induced from $V$
is denoted by $\mathrm{ind}_G(V)$.  Similarly, if $V$ is a
representation of $G$ the representation of $H$ formed by the
restriction of $V$ to $H$ is denoted by $\mathrm{res}_G(V)$. If $\chi$
is the character of $V$, then $\mathrm{ind}_G(\chi)$ denotes the
character of $\mathrm{ind}_G(V)$, and similarly for restrictions.
When it is clear which groups we are considering, they will not be
included in the notation. For a group $H \leq G$, the permutation
representation of $G$ on $G/H$ is the representation of $G$ induced by
the trivial representation of $H$ (see~\cite{MR1093239}, Section 1.12.)

For two characters $\phi$ and $\psi$ of a group $G$ there is an inner
product defined by 
\[ \langle \phi, \psi \rangle_G = \frac{1}{G} \sum_{g \in G} \phi(g)\psi(g).\] 
In particular, if $\pi$ is a permutation character of a group $G$, then 
\[\langle 1_G, \pi \rangle_G = \frac{1}{G} \sum_{g \in G} \fix(g).\]
This value is the number of orbits $G$ has on the set. 

Any character $\chi$ of $G$ can be expressed as a unique linear combination
of irreducible characters of $G$. If $\chi = \sum_{i=0}^k c_i
\chi_i$ where each $\chi_i$ is irreducible, then $\langle \chi ,\chi_i \rangle = c_i$.
A representation of $S_n$ is \textsl{multiplicity free} if no irreducible
character occurs more than once in its character decomposition.
If the permutation representation of $G$ in $S_n$ is multiplicity free, we say that $G$ is multiplicity free.
In particular, if $G$ is multiplicity free, then for every
partition $\la$ of $n$ $\langle \indg{G}{S_n},\chi_\la \rangle = 0$ or 1.

An important relationship with the inner product of induced and
restricted characters is \textsl{Frobenius reciprocity}.  Let $G$ and
$H$ be groups with $H \leq G$.  If $\phi$ is a character of $H$ and
$\psi$ is a character of $G$ then
\[
\langle \mathrm{ind}_G(\phi) , \psi \rangle_G 
= \langle \phi, \mathrm{res}_H(\psi) \rangle_H.
\]

\begin{theorem}\label{thm:basicfacts}
Let $n$ be a positive integer, and $H \leq G \leq S_n$. Then
\begin{enumerate}[(a)]
\item $\indg{H}{S_n} = \mathrm{ind}_{S_n}(\indg{H}{G})$;
\item The decomposition of $\indg{H}{G}$ is $1_G + \chi$ for some representation $\chi$ of $G$;
\item If the decomposition of $\indg{H}{S_n}$ is $\sum_{\la \vdash n} a_\la \chi^\la$, then 
      the decomposition of $\indg{G}{S_n}$ is $\sum_{\la \vdash n} b_\la \chi^\la$ where
      $b_\la \leq a_\la$;
\item If $H$ is multiplicity free then $G$ is multiplicity free. If $G$ is not multiplicity 
free then $H$ is not multiplicity free.\label{subgroups} \qed
\end{enumerate} 
\end{theorem}

\begin{proposition}\label{prop:intersectalt}
Let $G$ be a subgroup of $S_n$.  Assume that the decomposition of
$\indg{G}{S_n}$ into irreducible characters is $\sum_{\la
\vdash n} c_{\la} \chi^{\la}$.  If $G$ is not a subset of the alternating group
then 
\[
\indg{G\cap A_n}{S_n} = \sum_{\la \vdash n} c_{\la} (\chi^{\la} + \chi^{\la'}).
\]
\end{proposition} 
\proof
Let $C$ be a conjugacy class of $S_n$. Then the value of the character
of $\ind{G \cap A_n}$ at an element of $C$ is

\renewcommand{\arraystretch}{1.25} 
\begin{eqnarray*}
[S_n : (G \cap A_n)] \frac{|G \cap A_n \cap C|}{|C|} 
    &=& \left\{
     \begin{array}{ll} 
         2[S_n:G]\frac{|G \cap C|}{|C|}  & \textrm{if $C$ is even,}\\
         0 & \textrm{otherwise.} 
     \end{array}
       \right.
\end{eqnarray*}
\renewcommand{\arraystretch}{1} 

The value of the character of $\ind{G}$ on an element $g$ of $C$ is 
\[
[S_n:G]\frac{|G\cap C|}{|C|}.
\]
If $C$ is a class of even elements and $g$ is an element of $C$, then 
\[
(\chi^{\la}+\chi^{\la'})(g) = 2\chi^{\la}
\]
and if $C$ is a class of odd elements and $g \in C$, then 
\[
(\chi^{\la}+\chi^{\la'})(g) = 0.\qed
\]

If $G$ is a multiplicity-free subgroup of $S_n$, this result can be
used to determine when the subgroup $G \cap A_n$ is also multiplicity
free.  If the decomposition of $\indg{G}{S_n}$ is $\sum_{\la
\vdash n} c_\la \chi^\la$, then $G \cap A_n$ is multiplicity
free if and only if for every $\la \vdash n$, either $c_{\la}
= 0$ or $c_{\la'} = 0$.

\subsection{Littlewood-Richardson Rule}

Since we will make extensive use of the Littlewood-Richardson rule we
restate it here.  

First, let $V_1$ be a representation of $S_{d_1}$ and $V_2$ be a
representation of $S_{d_2}$, then $V_1 \otimes V_2$ is a
representation of $S_{d_1} \times S_{d_2}$.  So we can define $V_1
\circ V_2$ to be the representation of $S_{d_1 +d_2}$ induced from
$V_1 \otimes V_2$ of $S_{d_1} \times S_{d_2}$. 

The Littlewood-Richardson Rule gives the number of times an
irreducible representation (written as a partition of $n$) appears in
the decomposition of $V_1 \otimes V_2$ based on the irreducible
representations which occur in the decompositions of $V_1$ and $V_2$
(written as partitions of $d_1$ and $d_2$.)  To do this we need to
define for partitions $\nu$ and $\mu$ the \textsl{strict
$\nu$-expansion of $\mu$}.

Let $\mu$ and $\nu = [\nu_1,\nu_2,\dots ,\nu_k]$ be
partitions.  A \textsl{$\nu$-expansion} of $\mu$ is a partition that
is obtained in the following way. First, add $\nu_1$ boxes to $\mu$ in
such a way that no two boxes are in the same column and put a 1 in
each of these blocks. Next add $\nu_2$ blocks, again in such a way
that no two of these blocks are in the same column and place a 2 in
each of these blocks. Continue until $\nu_k$ blocks have been added
(these with a $k$ in each block).  This expansion is \textsl{strict}
if when the boxes are listed, going from right to left and starting at
the top and working down, for any $t \in \{1, \ldots , \sum_{i=1}^k
\nu_k\}$ each integer $p \in \{1,\ldots,k-1\}$ occurs at least as many
times as the integer $p+1$.

\begin{eg}
Below are all the $\mu$-expansions of $\la$ for $\mu
=[2,1]$ and $\la =[3]$. Only the last four are strict expansions.

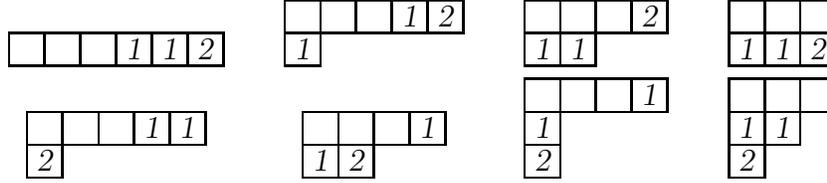
\begin{tabular}{cccc}
\begin{Young}
 & &  & 1 & 1 & 2\cr
\end{Young} &
\begin{Young}
 & & & 1 & 2\cr
1   \cr
\end{Young} &
\begin{Young}
 & & & 2 \cr 
1 & 1 \cr
\end{Young} &
\begin{Young}
 &  &  \cr 
1 & 1 & 2\cr
\end{Young} \\
\begin{Young}
 & &  & 1 & 1\cr 
2      \cr
\end{Young} &
\begin{Young}
 & & & 1 \cr
1  &  2     \cr
\end{Young} &
\begin{Young}
 & & & 1 \cr 
1 \cr
2 \cr
\end{Young} &
\begin{Young}
 &  &  \cr 
1 & 1 \cr
2 \cr 
\end{Young}
\end{tabular}
\end{eg}

The product $V_1 \circ V_2$ is associative and commutative, so
it is sufficient to state the decomposition for the case where $V_1$
and $V_2$ are irreducible representations. Assuming that $V_1$ is
irreducible, $V_1$ corresponds to a partition $\mu$ of $d_1$ and so we
will write $V_1$ as $V_\mu$. Similarly, we will write $V_2$ as $V_\nu$,
where $\nu$ is a partition of $d_2$.  Then
\[
V_\mu \circ V_\nu = \sum_{\la \vdash d_1 +d_2} N_{\mu,\nu}^\la V_\la.
\]
where the coefficients $N_{\mu,\nu}^\la$ are the number of ways $\mu$
can be expanded to $\la$ by a \textsl{strict $\nu$-expansion}.  A
proof of this rule can be found in~\cite{MR1354144} Chapter I, Section
9.

We can determine when the groups $S_k \times S_{n-k}$, $S_k \times
A_{n-k}$ and $A_k \times A_{n-k}$ are multiplicity-free subgroups of
$S_n$ with a straightforward application of the Littlewood-Richardson
Rule.

\begin{proposition}\label{prop:StimesS}
Let $k,n$ be integers with $2k \leq n$, then
\begin{enumerate}[(a)]
\item  $S_k \times S_{n-k}$ is multiplicity free,
\item  $A_k \times S_{n-k}$ is multiplicity free if and only if $k\neq 2$,
\item  $S_k \times A_{n-k}$ is multiplicity free if and only if $n-k\neq 2$,
\item  $A_k \times A_{n-k}$ is multiplicity free if and only if $k > 2$ and $2k \leq n-2$.
\end{enumerate}
\end{proposition}

\proof
The decomposition of $\indg{S_k}{S_k}$ into irreducible characters is
simply $[k]$. Similarly the decomposition of $\indg{S_{n-k}}{S_{n-k}}$
is $[n-k]$. All the $[n-k]$-expansions of $[k]$ are of the form $[k+i,
n-k-i]$ where $k+i \leq n-k-i$. By the Littlewood-Richardson Rule 
\[
\indg{S_k \times S_{n-k}}{S_n} = \sum_{i=0}^{k} [n-i,i].
\]

The decomposition of $\indg{A_k}{S_k}$ into irreducible
characters is $[k]+ [1^k]$. 
We can calculate the Littlewood-Richardson coefficients of $[k]
\circ [n-k]$ and $[1^k] \circ [n-k]$, we
have that 
\[
\indg{A_k \times S_{n-k}}{S_{n}} = [n-k,1^{k}]+[n-k+1,1^{k-1}] + \sum_{i=0}^{k} [n-i,i]
\]

Similarly 
\[
\indg{S_k \times A_{n-k}}{S_{n}} = [k,1^{n-k}]+[k+1,1^{n-k-1}] + \sum_{i=0}^{n-k} [n-i,i]
\]

Finally, by calculating the Littlewood-Richardson
coefficients, we have that 
\begin{align*}
\indg{A_k \times A_{n-k}}{S_n} = &[n-k,1^{k}]+  [n-k+1,1^{k-1}] +[k,1^{n-k}]\\ 
   &+[k+1,1^{n-k-1}] + \sum_{i=0}^{k} ([n-i,i]+[2^i,1^{n-2i}]). 
\end{align*}
If $k > 2$ and $2k < n-1$ then these characters are distinct.
\qed

In the following sections we give all multiplicity-free subgroups of
$S_n$. We first consider transitive subgroups, examining primitive and
imprimitive groups separately, and then the intransitive groups.

\section{Transitive Subgroups}
\label{sec:transitive}

Let $\Omega$ denote a set of $n$ elements and $G \leq S_n$ be a
permutation group on the set $\Omega$.  Let ${\Omega \choose k}$
denote the set of all unordered $k$-sets of distinct elements from
$\Omega$. Throughout this section $G$ will denote a subgroup of $S_n$
that is transitive on $\Omega$.

We will make extensive use of the following result from Saxl~\cite{MR627512}.  

\begin{lemma}\label{lem:saxl}
Let $k$ and $n$ be integers with $2(k+1) \leq n$.  If $G$ is a
multiplicity-free group with $\ell$ orbits on ${\Omega \choose k}$,
then $G$ has either $\ell$ or $\ell+1$ orbits on ${\Omega \choose
k+1}$.
\end{lemma} 

\proof
Let $\la = (\la_1,\la_2,\dots ,\la_k)$ be a
partition of $n$ and let 
\[
S_\la = S_{\la_1} \times S_{\la_2} \times  \cdots \times S_{\la_k}
\]
be the Young subgroup corresponding to $\la$.

By Frobenius reciprocity
\[
\langle \ind{G}, \ind{S_\la} \rangle_{S_n}
= \langle 1, \mathrm{res}(\indg{S_\la}{S_n}) \rangle_G.
\]
The value of the inner product $\langle 1, \mathrm{res}(\indg{S_\la}{S_n})
\rangle_G$ is the number of orbits of the action of $G$ on
the cosets $S_n / S_\la$.
For the partition $\la =[n-k,k]$, the value of $\langle 1, \mathrm{res}(\indg{S_\la}{S_n})
\rangle_G$ is the number of $G$ orbits on ${\Omega \choose k}$.
Since the group $G$ has exactly $\ell$ orbits on ${\Omega \choose k}$
it follows that
\[\langle 1, \mathrm{res}(\ind{S_{n-k} \times S_k}) \rangle_G = \ell.\]

Since $2(k+1) \leq n$, from the proof of Proposition~\ref{prop:StimesS}
\[
\indg{ S_{n-k-1}\times S_{k+1} }{S_n} = \indg{S_{n-k} \times S_k}{S_n} + [n-k-1,k+1].  
\]
The number of orbits $G$ has on ${\Omega \choose k+1}$
is 
\begin{align*}
\langle 1, \mathrm{res}(\ind{S_{n-(k+1)}\times S_{k+1}}) \rangle_G 
&= \langle 1, \mathrm{res}(\ind{S_{n-k}\times S_k }) \rangle_G + \langle 1, \mathrm{res}(\chi^{[n-k,k]}) \rangle_G \\
&= \ell + \langle 1, \mathrm{res}(\chi^{[n-k,k]}) \rangle_G \\
&= \ell + \langle \ind{G}, \chi^{[n-k,k]} \rangle_{S_n}.
\end{align*}

Since $G$ is multiplicity free, $\langle \indg{G}{S_n}, \chi^{[n-k,k]}
\rangle=0$ or $1$ and $G$ has either $\ell$ or $\ell+1$ orbits on
${\Omega
\choose k+1}$.
\qed

The next result clearly follows from Lemma~\ref{lem:saxl}.
\begin{corollary}
Let $n$ and $k$ be integers.
If $G \leq S_n$ is multiplicity free, then the number of orbits of $G$ on
${\Omega \choose k}$ is no more than $k+1$. Further, if $G$ is
transitive and multiplicity free, then the number of orbits of $G$ on
${\Omega \choose k}$ is no more than $k$. \hfill \qed
\end{corollary}

\subsection{Primitive Subgroups}
\label{sec:primitive}

In this section we consider the primitive transitive subgroups of the
symmetric group. In the following section we will consider the
imprimitive transitive subgroups of the symmetric group.

For $n\geq 6$, if a multiplicity-free group $G \leq S_n$ is primitive
and transitive then it is 2-transitive (\cite{MR875503}, Lemma 4.1.)
Saxl~\cite{MR627512} shows that no 2-transitive subgroup with degree
18 or greater is multiplicity free.

We only consider primitive subgroups with degree less than or equal to
18.  A list of all primitive groups of degree less than or equal to 20
is given by Sims in Table 1 of~\cite{MR0257203}. 

With GAP, for each of these groups it is possible to construct the
representation induced on $S_n$ by the trivial character on the
group. Then, using GAP, we obtain a list of the irreducible characters
on $S_n$ which occur in the induced character. The multiplicity of
each irreducible character in the induced character can be found by
taking the inner product of the irreducible character with the induced
character. If all the multiplicities are equal to 1 then the induced
character is multiplicity free. Table~\ref{table:transitive} is a
complete list of the primitive, multiplicity-free groups with degree
less than 18.

\subsection{Imprimitive Subgroups}
\label{sec:imprimitive}
In this section we consider imprimitive, transitive subgroups of the
symmetric group.  Since $G$ is transitive, it has one orbit on
$\Omega$ and at most $k$ orbits on ${\Omega \choose k}$ for all $k
\geq 1$.

\begin{lemma}\label{lem:blockstructure}
If $G \leq S_n$ is a multiplicity-free group that is imprimitive, then
the blocks of imprimitivity have one of the following structures:
\begin{enumerate}[(a)]
\item{2 blocks of size $n/2$;}
\item{$n/2$ blocks of size 2;}
\item{3 blocks of size no more than 5;}
\item{at most 5 blocks of size 3.}
\end{enumerate}
\end{lemma}

\proof
It is sufficient to show that if the blocks of imprimivity do not have
one of the structures above, then the group $G$ is not multiplicity free. 

Assume that $G$ has three blocks of imprimitivity.  For any $k$-set of
elements from $\Omega$, let the unordered triple $(x,y,z)$ denote the
number of elements in the $k$-set from three distinct, but unordered,
blocks of imprimitivity. Distinct triples $(x,y,z)$ correspond to
distinct orbits of $G$ on ${\Omega \choose k}$. If the blocks of
imprimitivity have at least six elements, then there are seven
distinct unordered triples $(x,y,z)$ and the group $G$ has at least
seven orbits on ${\Omega \choose 6}$. By Lemma~\ref{lem:saxl} this
means that $G$ is not multiplicity free.

Next assume that each of the blocks of imprimitivity has size three
and that there are six such blocks.  By counting the unordered
triples $(x,y,z)$ as above, $G$ has seven orbits on ${\Omega \choose
6}$ and by Lemma~\ref{lem:saxl} the group $G$ is not multiplicity free.

Finally, if $G$ has at least four blocks of imprimitivity each of size
at least four, then there are at least five orbits on ${\Omega \choose
4}$.  Again, by Lemma~\ref{lem:saxl} the group $G$ is not multiplicity free.\qed

The next result shows that an imprimitive group is a subgroup of a
wreath product. Let $G$ be a group and $k$ a positive integer. The
direct product of $k$ copies of $G$ is denoted $G^k$. The wreath product
of $G$ and $S_k$, denoted $G \wr S_k$, is the group with elements from $G^k
\times S_k$ with multiplication defined by
\[
 (g_1,g_2,\ldots , g_k\! :\! \pi)(h_1,h_2,\ldots , h_k\!:\! \tau)
 :=(g_1 h_{\pi^{-1}(1)}, g_2 h_{\pi^{-1}(2)},\ldots, g_k h_{\pi^{-1}(k)}\!:\!\pi\tau).
\]

\begin{corollary}
If $G$ is a transitive, imprimitive and multiplicity-free subgroup of $S_n$
then $G$ is a subgroup of $S_\ell \wr S_k$ where $(\ell,k)$ is one of the 
following pairs:
\begin{enumerate}[(a)]
\item $(\ell,k) = (2,k)$;
\item $(\ell,k) = (\ell,2)$;
\item $(\ell,k)$ is one of $(3,3),(4,3)$ or $(5,3)$;
\item $(\ell,k) = (3,4)$. 
\end{enumerate}
\end{corollary}

\proof
Let $B=\{B_1,B_2,\dots , B_k\}$ be the blocks of imprimitivity for $G$
with $|B_i|=\ell$ for all $i\in \{1,\dots, k\}$.  Set $H_i =
\stab{B_i}{G}$ to be the setwise stabliser of $B_i$ in $G$ and $K$ to
be the permutation group induced by the action of $G$ on the blocks of
imprimitivity.  Since $G$ is transitive, $H_i =H_j$ for all $i,j \in
\{1,\dots, k\}$, so we can write $H$ for $H_i$.

Then $H_{\upharpoonright B_i} \leq S_\ell$ and $K \leq S_k$ and
$H_{\upharpoonright B_i} \wr K \leq S_\ell \wr S_k$. The group $G$ can
be embedded in $H_{\upharpoonright B_i} \wr K$ and $G \leq S_\ell \wr
S_k$.

Finally, if $G$ has three blocks each of size five, then $G$ is a
subgroup of $S_5 \wr S_3$. This group is not multiplicity free (this
can be checked using GAP) and $G$ is not multiplicity free by
Theorem~\ref{thm:basicfacts}(d).
\qed

It is possible to further restrict which imprimitive groups are
multiplicity free by extending the argument used in
Lemma~\ref{lem:blockstructure}. For integers $t \leq n$, a group
acting on an $n$-set is \textsl{$t$-homogeneous} if it has one orbit
on the collection of unordered $t$-sets from the $n$-set.

\begin{corollary}\label{cor:thomonblocks}
Let $G$ be an imprimitive multiplicity-free group and let
$B=\{B_1,B_2,\dots , B_k\}$ be the blocks of imprimitivity. Assume
each block has size $\ell$. Let $K$ be the permutation group induced
by the action of $G$ on the blocks of imprimitivity. Then 
\begin{enumerate}[(a)]
\item{the group $K$ is $t$-homogeneous on the blocks for all $t \leq k$;}
\item{the group $\stab{B_i}{G}$ of $G$ restricted to the block $B_i$ is
$t$-homogeneous for all $t \leq \ell$;}
\item {$\stab{B_i}{G} = \stab{B_j}{G}$ for all $i,j \in \{1, \dots, k\}$.}
\end{enumerate}
\end{corollary}
\proof
The first two results can be seen by
counting the number of orbits of $G$ on ${\Omega \choose 3}$.

Since $G$ is transitive there is a $\alpha \in K$
exchanges the blocks $B_i$ and $B_j$. If $g \in \stab{B_i}{G}$ 
then $\alpha g \alpha^{-1} \in \stab{B_j}{G}$.
\qed

Beaumont and Peterson~\cite{MR0067889} have determined the
subgroups of $S_n$ that are $t$-homogeneous for all $t \leq n$.  They
are $S_n$, $A_n$, $AGL(1,5)$ with $n=5$, $PGL(2,5)$ with $n=6$, and
$PGL(2,8)$ or $P\Gamma L(2,8)$ with $n=9$.

Next we will determine, for each block structure, which imprimitive
groups are multiplicity free.

\subsubsection{Multiplicity-free subgroups of $S_k \wr S_2$}

In this section, we give all multiplicity-free subgroups of $S_k \wr
S_2$, starting with the group $S_2 \wr S_2$. 

\begin{proposition}\label{prop:subgroupsS2wrS2}
The group $S_2 \wr S_2$ is multiplicity-free and
the only transitive, proper subgroup of $S_2 \wr S_2$ 
which multiplicity-free in $S_4$ is the group generated by the permutation $(1,3,2,4)$.
\end{proposition}
\proof
Using GAP, it is straightforward to calculate that
\[
\indg{S_2 \wr S_2}{S_4} =[4]+[2,2].
\]
Further, it is not difficult to find all the multiplicity-free
subgroups of $S_4$, these are:
\[
A_3,\ S_3,\ S_2 \times S_2,\ \langle(1,3,2,4) \rangle,\ S_2 \wr S_2,\ A_4,\ S_4.
\]
Of these, only $\langle(1,3,2,4) \rangle$ is a proper, transitive subgroup of 
$S_2\wr S_2$.\qed

The subgroup of $S_2 \wr S_2$ generated by the permutation $(1,3,2,4)$
is part of a family of subgroups of the wreath product that is
multiplicity free in several cases. For a group $G$, recall that $G^k$ denotes
the direct product of $k$ copies of $G$.  Each element of $\sigma \in
S_k$ defines an automorphism of $G^k$ by $\sigma(g_1,g_2, \cdots ,
g_k) = (g_{\sigma^{-1}(1)},g_{\sigma^{-1}(2)},\cdots
,g_{\sigma^{-1}(k)})$, and thus there is a homomorphism of $S_k$ to
the automorphism group of $G^k$.  With this homomorphism, the wreath
product $G \wr S_{k}$ is a semi-direct product, $G^k \rtimes S_k$.
For any subgroup $H$ of $G^k$ which is normal in $G\wr S_k$, the
semidirect product of $S_k$ and $H$, using the above homomorphism,
gives a subgroup of $G \wr S_k$. We will denote this group by $H
\rtimes S_k$. For example, the subgroup generated by the permutation
$(1,3,2,4)$ can also be expressed as the more complicated $(S_2^2 \cap
A_4) \rtimes S_2$.

The group $S_\ell \wr S_2$ is multiplicity free and its decomposition
is known. (See \cite{MR627512} Example 2.3.)

\begin{proposition}\label{prop:SkwrS2}
The group $S_\ell \wr S_2$ is multiplicity free in $S_{2\ell}$ with
\[
\indg{S_\ell \wr S_2}{S_{2\ell}} = \sum_{i=0}^{\lfloor \ell/2 \rfloor } [2\ell-2i,2i].\qed
\] 
\end{proposition} 

Next we will determine all the multiplicity-free subgroups of $S_\ell
\wr S_2$. To do this we will need two results. First is the list of subgroups of $S_2 \wr
S_2$ whose permutation representations are multiplicity free in $S_2
\wr S_2$ and second is the list of characters of $S_\ell \wr S_2$
induced from the subgroups of index 2. The subgroups which are
multiplicity free in $S_2 \wr S_2$ are given in the following poset
ordered by inclusion.

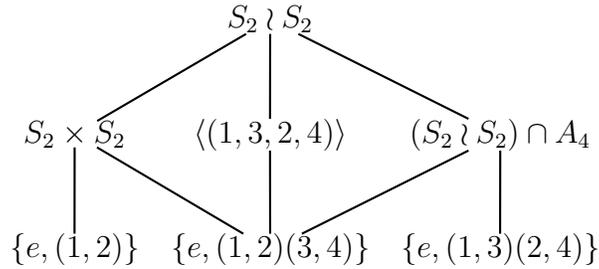
\begin{center}
\begin{tabular}{ccc}
 & \rnode{a}{$S_2 \wr S_2$} &   \\
 & & \\ & & \\
\rnode{b}{$S_2 \times S_2$} & \rnode{c}{$\langle(1,3,2,4) \rangle$} 
  & \rnode{d}{$(S_2 \wr S_2 ) \cap A_4$}\\
 & & \\ & & \\
\rnode{e}{$\{e,(1,2)\}$} & \rnode{f}{$\{e, (1,2)(3,4)\}$} & \rnode{g}{$\{e,(1,3)(2,4)\}$} \\
\end{tabular}
\end{center}
\ncline{-}{a}{b}\ncline{-}{a}{c}\ncline{-}{a}{d}
\ncline{-}{b}{f}\ncline{-}{c}{f}\ncline{-}{d}{f}
\ncline{-}{b}{e}\ncline{-}{d}{g}

\begin{proposition}\label{prop:linearcharsofwreath}
Let $\ell$ be a positive integer. 
\begin{enumerate}[(a)]
\item The decomposition of $\indg{ (S_\ell \wr S_2) \cap A_{2\ell}}{S_\ell
\wr S_2}$ is $1 + \sigma$ where $\sigma$ is the sign character. 
\item The decomposition of $\indg{ S_\ell \times S_\ell }{S_\ell \wr S_2}$
is $1 + \psi$ where $\psi$ is a linear character. 
\item The decomposition of $\indg{ ((S_\ell \times S_\ell) \cap A_{2\ell})
\rtimes S_2}{S_\ell \wr S_2}$ is $1 + \sigma\psi$. 
\end{enumerate}
\end{proposition}
\proof 
Part (a) follows from Proposition~\ref{prop:intersectalt} and part (b)
follows from the fact that the index of $S_\ell \times S_\ell$ in
$S_\ell \wr S_2$ is two.

To prove Part (c), we will prove that the group $S_\ell \wr S_2$ has
exactly four irreducible linear characters, 1, $\sigma$, $\psi$ and
$\sigma \psi$.  

Each linear character of $S_\ell \wr S_2$ is the character of the
permutation representation of $S_\ell \wr S_2$ on the cosets $(S_\ell
\wr S_2)/N$ for a normal subgroup $N$.  If $N$ is normal in $S_\ell
\wr S_2$, then $A_\ell \times A_\ell \leq N$.  The group $(S_\ell \wr S_2)/(A_\ell \times
A_\ell)$ is isomorphic to $S_2 \wr S_2$, which is the dihedral group
of order 8. Since the dihedral group of order eight has four linear
characters, $N/(A_\ell
\times A_\ell)$ is the kernel of one of these four characters.  Thus
there are at most four subgroups $N$ and $S_\ell \wr S_2$ has exactly
four linear characters.
\qed 

\begin{theorem}\label{thm:subgroupsSkwrS2}
For $\ell > 2$ and $\ell \neq 5,6$ and $9$ the transitive, multiplicity-free, proper
subgroups of $S_\ell \wr S_2$ with their decompositions are exactly the following groups:
\begin{enumerate}[(a)]
\item $(S_\ell \wr S_2) \cap A_{2\ell}$,
\item $((S_\ell \times S_\ell) \cap A_{2\ell})\rtimes S_2$,
\item $A_\ell \wr S_2$.
\end{enumerate}
\end{theorem}
\proof
Assume $G$ is a multiplicity-free subgroup of $S_\ell \wr S_2$.
Since the stabiliser $\stab{B_i}{G}$ of $G$ restricted to a
block $B_i$ is $t$-homogeneous for all $t \leq \ell$ and for the given
values of $\ell$ this group must be either $A_\ell$ or $S_\ell$.
In particular, for the given values of $\ell$, $A_\ell \times A_\ell \leq G$

The group $A_\ell \times A_\ell$ is a normal subgroup of $G$ so we can
define the group $\overline{G} = G/(A_\ell \times A_\ell)$. Then
$\overline{G}$ is a subgroup of $(S_\ell \wr S_2) / (A_\ell \times
A_\ell) \cong S_2 \wr S_2$.

Since $G$ is multiplicity-free in $S_{2\ell}$, the representation $\indg{G}{S_\ell
\wr S_2}$ is also multiplicity free (Theorem~\ref{thm:basicfacts}(a)). This implies that
$\indg{\overline{G}}{(S_\ell \wr S_2)/(A_\ell \times A_\ell)}$ is a
multiplicity-free representation of $S_2 \wr S_2$. The proper
multiplicity-free subgroups of $S_2 \wr S_2$ are given in the proof of
Proposition~\ref{prop:subgroupsS2wrS2}.

This means that $G$ is a transitive group with the property that
$\overline{G}$ is a multiplicity-free subgroup of $S_2 \wr S_2$.
There are only three such groups
\[
(S_\ell \wr S_2) \cap A_{2\ell} \quad ((S_\ell \times S_\ell) \cap A_{2\ell})\rtimes S_2
\quad A_\ell \wr S_2. 
\]
We need to prove that these three groups are multiplicity free.

By Proposition~\ref{prop:intersectalt} and
Proposition~\ref{prop:SkwrS2},
\[
\indg{(S_\ell \wr S_2) \cap A_{2\ell}}{S_{2\ell}} = \sum_{i=0}^{\lfloor \ell/2 \rfloor } 
 \left( [2\ell-2i,2i] +  [2^{2i},1^{2\ell-4i}]\right).
\]
For $\ell>2$, the group $(S_\ell\wr S_2) \cap A_{2\ell}$ is multiplicity free. 
Further, since the decomposition of $\indg{ (S_\ell \wr S_2) \cap A_{2\ell}}{S_\ell
\wr S_2}$ is $1 + \sigma$, we also have that 
\[
\mathrm{ind}_{S_{2\ell}}( \sigma) = \sum_{i=0}^{\lfloor \ell/2 \rfloor } [2^{2i},1^{2\ell-4i}].
\]

The decomposition of $\indg{ S_\ell \times S_\ell }{S_\ell \wr S_2}$
is $1 + \psi$ where $\psi$ is a linear character.  In
Proposition~\ref{prop:StimesS} the decomposition of $\indg{ S_\ell
\times S_\ell }{S_{2\ell}}$ is given. Using this, together with the
decomposition of $\indg{S_\ell \wr S_2}{S_{2\ell}}$ from
Proposition~\ref{prop:SkwrS2}, we have that
\[\mathrm{ind}_{S_{2\ell}}(\psi) = \sum_{i=0}^{\lfloor \ell/2 \rfloor}[2\ell-2i-1, 2i+1].\]

The decomposition of $\indg{ ((S_\ell \times S_\ell) \cap A_{2\ell}) \rtimes S_2}{S_\ell
\wr S_2}$ is $1 + \sigma\psi$.  Since $\sigma$ is the sign character,
\[\mathrm{ind}_{S_{2\ell}}(\sigma\psi) = 
   \sum_{i=0}^{\lfloor \ell/2 \rfloor} [2^{2i+1},1^{2\ell-4i-2}],\]  
from this it follows that
\[\indg{ ((S_\ell \times S_\ell) \cap A_{2\ell}) \rtimes S_2}{S_\ell \wr
 S_2}= \sum_{i=0}^{\lfloor \ell/2 \rfloor}\left( [2\ell-2i, 2i]+[2^{2i+1},1^{2\ell-4i-2}]\right).\]  
For all $\ell$ the group $((S_\ell \times S_\ell) \cap
A_{2\ell}) \rtimes S_2$ is multiplicity free.

Since $A_\ell \times S_\ell$ is a subgroup of $S_\ell \times S_\ell$,
we see that
\[
\indg{ A_\ell \times S_\ell}{S_\ell \wr S_2} =1 + \psi + \phi
\]
where $\phi$ is a 2 dimensional character.  {}From the
Littlewood-Richardson Rule (Proposition~\ref{prop:StimesS}), the
decomposition of $A_\ell \times S_\ell$ is known and from this we can
deduce that 
\[
\mathrm{ind}_{S_{2\ell}}(\phi) = [\ell+1, 1^{\ell-1}]+[\ell,1^\ell].
\]

Similarly, we know the decomposition of $A_{\ell} \times A_\ell$.
From this, and the fact that $A_\ell \times A_\ell$ is a subgroup of
the three groups 
\[
S_\ell \times S_\ell,\quad
	(S_\ell \wr S_2) \cap A_{2\ell},\quad
	(S_\ell \times S_\ell)\cap A_{2\ell})\rtimes S_2,
\]
we conclude that 
\[
\indg{A_\ell \times A_\ell}{S_\ell \wr S_2} =1 + \sigma + \psi +\sigma \psi + 2\phi.
\]

Finally, $A_\ell \times A_\ell$ is a subgroup of $A_\ell \wr S_2$ so
the decomposition of $\indg{A_\ell \wr S_2}{S_{2\ell}}$ is
a subset of the characters in $1 + \sigma + \psi + \sigma \psi + 2\phi$.  Since $A_\ell
\wr S_2$ is transitive, the decomposition of $\indg{A_\ell \wr
S_2}{S_{2\ell}}$ does not include the partition $[2\ell -1,1]$, this
means that $\psi$ is not in the decomposition. Since 
$[S_\ell \wr S_2: A_\ell \wr S_2]= 4$ and 
$\phi$ has degree 2 while $\sigma$ has
degree 1, we conclude that 
\[
\indg{A_\ell \wr S_2}{S_\ell\wr S_2} =1 + \sigma\psi + \phi.
\]
The decomposition of each of these characters is known, and
\[
\indg{A_\ell \wr S_2}{S_{2\ell}} = \sum_{i=0}^{\lfloor \ell/2 \rfloor } 
 \left( [2\ell-2i,2i] +[2^{2i+1},1^{2\ell-4i-2}] \right) +
 [\ell+1, 1^{\ell-1}]+[\ell,1^\ell].
\]
It is clear from the formula that $A_\ell \wr S_2$ is multiplicity free for all $\ell>2$.
\qed

It is interesting to note that each of the representations induced from the linear characters $1$,
$\sigma$, $\psi$ and $\sigma\psi$ of $S_\ell \wr S_2$ are multiplicity
free in $S_{2\ell}$.

\begin{theorem}\label{thm:specialwreaths}
For $\ell = 5,6$ or 9 the transitive, multiplicity-free,
subgroups of $S_\ell \wr S_2$ are the subgroups given in Theorem~\ref{thm:subgroupsSkwrS2}
and the following:
\begin{enumerate}[(a)]
\item{$AGL(1,5) \wr S_2$ for $\ell =5$,}
\item{$PGL(2,5) \wr S_2$ for $\ell =6$,}
\end{enumerate}
\end{theorem}

\proof
Assume $G \leq S_{2\ell}$ is multiplicity free.  If $\stab{B_i}{G}
=S_n$ or $A_n$, then $G$ is one of the subgroups given in
Theorem~\ref{thm:subgroupsSkwrS2}.

{}From Corollary~\ref{cor:thomonblocks} there are four cases to consider,
first where $\ell =5$ and $\stab{B_i}{G} = AGL(1,5)$,
second where $\ell =6$ and $\stab{B_i}{G} = PGL(1,5)$ and finally 
 where $\ell =9$ and $\stab{B_i}{G} = PGL(2,8)$ or $P\Gamma L(2,8)$.

Assume $\ell =5$ and $\stab{B_i}{G} = AGL(1,5)$.  This means 
\[
AGL(1,5) \times AGL(1,5) \leq G.
\]
If $G$ is multiplicity free in $S_{10}$, then
$G$ must also be multiplicity free in $S_5 \wr S_2$.  Further, the
group $G/(AGL(1,5) \times AGL(1,5))$ must also be multiplicity free in
$(S_5 \wr S_2) / (AGL(1,5) \times AGL(1,5))$.  This last group is
congruent to $S_3 \wr S_2$ and from Proposition~\ref{prop:SkwrS2} we
have a complete list of its multiplicity-free subgroups.  From this
list, we can build the list of the groups $G$, for which $G/(AGL(1,5)
\times AGL(1,5))$ is one of the multiplicity-free subgroups of $S_3
\wr S_2$.  Using GAP, it is possible to check which of these groups
are multiplicity free in $S_{10}$.

For $\ell =6$ the proof is similar since 
\[(S_6 \wr S_2) / (PGL(2,5) \times PGL(2,5)) \cong S_3 \wr S_2.\]

The group $P\Gamma L(2,8)$ is a maximal subgroup of $A_9$.  If $G$ is
transitive and $PGL(2,8) \times PGL(2,8) \leq G$ then $G \leq P\Gamma
L(2,8) \wr S_2$.  Since the group $P\Gamma L(2,8) \wr S_2$ is not
multiplicity free, no subgroups of it are multiplicity free
(Theorem~\ref{thm:basicfacts}(\ref{subgroups}).)
\qed

\subsubsection{Multiplicity free subgroups of $S_2 \wr S_k$}

In this section we find all multiplicity-free subgroups of $S_2 \wr
S_k$.  We start by stating the decomposition of $S_2 \wr S_k$, for a
proof see \cite{MR627512}, Example 2.2.

\begin{proposition}\label{prop:wreath1}
The group $S_2 \wr S_k$ is multiplicity free in $S_{2k}$.  Moreover,
for a partition $\la$ of $k$ with $\la=[\la_1,\la_2,
\dots, \la_a]$, define a partition of $2k$ by 
$2\la = [2\la_1, 2\la_2, \dots, 2\la_a]$, then
\[
\indg{S_2 \wr S_k}{S_{2k}} = \sum_{\la \vdash k} 2\la. \qed
\]
\end{proposition} 

Let $G$ be a subgroup of $S_2\wr S_k$ with $B=\{B_1,B_2,\dots,B_k\}$
its blocks of imprimitivity. Let $K$ be the permutation group induced
by the action of $G$ on the blocks of imprimitivity. If $G$ is
multiplicity free, then by Corollary~\ref{cor:thomonblocks},
$\stab{B_i}{G} =S_2$ and $K$ is $t$-homogeneous
for all $t \leq k$.  The group $K$ is $S_n$ or
$A_n$ or one of the following groups: $AGL(1,5)$, $PGL(2,5)$,
$(PGL(2,8)$ or $P \Gamma L(2,8)$~\cite{MR0067889}. With GAP we can
check that the wreath product of $S_2$ with each of the last four
groups is not multiplicity free. If $K$ is one of
the last four groups, then $G$ is not multiplicity free since it is a
subgroup of a group that is not multiplicity free.

Thus $A_k \leq K$ for all $k$ and $G$ has index at most 4 in 
$S_2 \wr S_k$~\cite{MR627512}. There are only four proper
subgroups of $S_2 \wr S_k$ which could be multiplicity free, these are:
\[
(S_2\wr S_k) \cap A_{2k} \quad ( (S_2)^k  \cap A_{2k})\rtimes S_k 
\quad S_2\wr A_k \quad ( (S_2)^k  \cap A_{2k})\rtimes A_k.
\]

\begin{proposition}
The group $(S_2\wr S_k) \cap A_{2k}$ is multiplicity free if and only if $k$ is odd.
\end{proposition}
\proof
If $k$ is even, the conjugate partitions $[2^k]$ and $[k,k]$ both
occur in the decomposition of $\ind{S_2\wr S_k}$. By Proposition~\ref{prop:intersectalt},
$(S_2\wr S_k) \cap A_{2k}$ is not multiplicity free.

If $(S_2\wr S_k) \cap A_{2k}$ is not
multiplicity free then there must be a pair of conjugate characters in
the decomposition of $S_2\wr S_k$. By
Proposition~\ref{prop:wreath1}, there must be two partitions 
$\la=[\la_1,\la_2,\dots,\la_a]$
and $\mu = [\mu_1,\mu_2,\dots,\mu_b]$ of $k$ with 
\[
(2 \la)' = [2\la_1,2\la_2,\dots,2\la_a]' =[2\mu_1,2\mu_2,\dots,2\mu_b] =2\mu.
\]  
By definition of conjugation, $a=2\mu_1$, and $a$ must be even.
Moreover, if $\la_{i-1} > \la_{i}$ then $i-1$ must occur in $2\mu$ and hence $i$ must be odd.
Therefore, if $i$ is even then $\la _{i-1} = \la_{i}$.
Thus 
\begin{eqnarray*}
k &=& \sum_{i=1}^{a} \la_i 
 = \sum_{\substack{i=2 \\ i \hspace{1mm}\mathrm{even}}}^{a} (\la_{i-1} + \la_i)
 = \sum_{\substack{i=2 \\ i \hspace{1mm}\mathrm{even}}}^{a} 2 \la_i
\end{eqnarray*}
and $k$ is even.
\qed

Next we determine all $k$ for which the group $S_2 \wr A_k$ is
multiplicity free.  The decomposition of $\indg{S_2 \wr A_k}{S_2 \wr
S_k}$ is $1+\psi_k$ where $\psi_k$ is the sign character on the blocks
of imprimitivity of $S_2 \wr S_k$. We will give the decomposition of
$\mathrm{ind}_{S_{2k}}(\psi_k)$.

First we need to introduce some notation.  Let $\la$ be a partition of
an integer $n$. For a node $(i,i)$ in the Young diagram of $\la$ the
\textsl{diagonal hook}, $h_{i}$, is the set of all nodes $\{(i,j) :
j\geq i\} \cup \{(k,i) : k \geq i\}$. The \textsl{width}, $W(h_i)$ of
a diagonal hook $h_i$ is $|\{(i,j) : j\geq i\}|$ and the
\textsl{depth}, $D(h_i)$ is $|\{(k,i) : k\geq i\}|$.

We will also need the following result which was observed by Saxl
in~\cite{MR0450385} as a direct result of Mackey's subgroup theorem
(see \cite{MR0144979}, Section 10.13 ). 

\begin{lemma}\label{lem:mackey}
If $G$ is a permutation group on a set $\Omega$ and $H$ is a transitive subgroup,
then for $a \in \Omega$ and $\chi$ an irreducible character of $\stab{G}{a}$ 
\[  \mathrm{res}_H(\mathrm{ind}_G(\chi)) = \mathrm{ind}_H (\mathrm{res}_{\stab{H}{a}}(\chi)).\]
\end{lemma}

Before stating the the decomposition of
$\mathrm{ind}_{S_{2k}}(\psi_k)$ we give an example.

\begin{eg}
The decomposition of $\mathrm{ind}_{S_{12}}(\psi_6)$ is
$[4,4,4]+[5,4,2,1]+[6,3,1,1,1]+[7,1,1,1,1,1]$. Below are the Young
tableaux for these partitions. In  each tableaux the diagonal hooks
marked with the symbols, $\circ$, $\ast$ and $\cdot$.

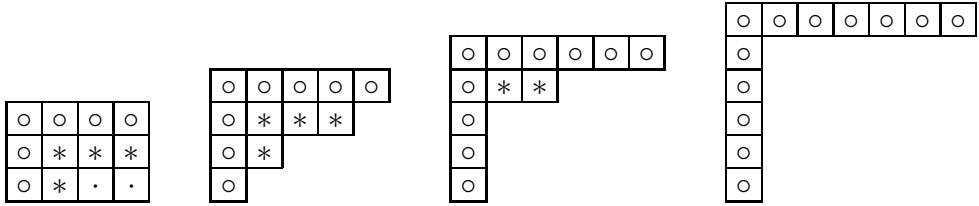
\begin{center}
\begin{tabular}{cccc}
\begin{Young}
$\circ$ & $\circ$&$\circ$ &$\circ$ \cr
$\circ$ & $\ast$ &$\ast$ &$\ast$  \cr
$\circ$ & $\ast$ & $\cdot$ & $\cdot$ \cr
\end{Young} &
\begin{Young}
$\circ$ &$\circ$ &$\circ$ &$\circ$ & $\circ$ \cr
$\circ$ & $\ast$&$\ast$ & $\ast$ \cr
$\circ$ & $\ast$ \cr
$\circ$ \cr
\end{Young} &
\begin{Young}
$\circ$ &$\circ$ &$\circ$ &$\circ$ &$\circ$ & $\circ$ \cr
$\circ$ &$\ast$ & $\ast$ \cr
$\circ$ \cr$\circ$ \cr$\circ$ \cr
\end{Young} &
\begin{Young}
$\circ$ &$\circ$ &$\circ$ &$\circ$ &$\circ$ &$\circ$ & $\circ$ \cr
$\circ$\cr$\circ$ \cr$\circ$\cr$\circ$\cr $\circ$\cr
\end{Young} 
\end{tabular}
\end{center}
\end{eg}

\begin{theorem}\label{thm:S2wrAk}
The decomposition of $\mathrm{ind}_{S_{2k}}(\psi_k)$ is the sum of all
partitions of $2k$ with the property that every diagonal hook in the
partition has depth one less than its width.
\end{theorem}
\proof
The following proof uses the method in Example 2.1 in~\cite{MR627512}
(this result is our Proposition~\ref{prop:wreath1}.)
In \cite{MR627512}, this method is attributed to James and Saxl.

We prove this theorem by induction on $k$. First, if $k=3$ then
decomposition of 
\[
\mathrm{ind}_{S_{6}}(\psi_3) =[4,1,1]+[3,3].
\]
A simple check of the Young diagrams for all partitions of 6 shows that
the theorem holds in this case.

First we give a description of the set of partitions in the
decomposition of $\mathrm{res}_{S_{2k-1}}(\mathrm{ind}_{S_{2k}}(\psi_k))$.  In the
second part of the proof, we show that if there is partition in
decomposition of $\mathrm{ind}_{S_{2k}}(\psi_k)$ that is not of the
right form, then the decomposition of the restriction of the corresponding
character to $S_{2k-1}$ is not contained in this set.

Apply Lemma~\ref{lem:mackey} with $\Omega$ the set of uniform
$k$-partitions of a $2k$-set, $G=S_{2k}$ and $H=S_{2k-1}$. Then
$\stab{G}{a} = S_2 \wr S_k$ and $\stab{H}{a} = S_2 \wr S_{k-1}$ for
any $a \in \Omega$ and thus
\begin{eqnarray*}
\mathrm{res}_{S_{2k-1}}(\mathrm{ind}_{S_{2k}}(\psi_k))
&=& \mathrm{ind}_{S_{2k-1}}( \mathrm{res}_{S_2 \wr S_{2k-1}} \psi_k).
\end{eqnarray*}

{}From the definition of $\psi_k$, we have that
\begin{align*}
\mathrm{ind}_{S_{2k-1}}(\mathrm{res}_{S_2 \wr S_{k-1}}{\psi_k})
&=  \mathrm{ind}_{S_{2k-1}}(\psi_{k-1}) \\
&=  \mathrm{ind}_{S_{2k-1}} ( \mathrm{ind}_{S_{2k-2}}({\psi_{k-1}})). 
\end{align*}
These together imply that
\begin{eqnarray}\label{eqn:equivdecomps}
\mathrm{res}_{S_{2k-1}}(\mathrm{ind}_{S_{2k}}(\psi_k))
=\mathrm{ind}_{S_{2k-1}} ( \mathrm{ind}_{S_{2k-2}}({\psi_{k-1}})).
\end{eqnarray}

By the induction hypothesis, $\mathrm{ind}_{S_{2k-2}}({\psi_{k-1}})$
is the sum of all partitions of $2k-2$ with the property that every
diagonal hook in the partition has depth one less than its width.

Using the Littlewood-Richardson rule, we see that the decomposition of
$\mathrm{ind}_{S_{2k-1}}(\mathrm{ind}_{S_{2k-2}}({\psi_{k-1}}))$ is
the sum of all partitions of $2k-1$ with all diagonal hooks, except one, with
depth one less then its width.  This one hook will have depth equal to
width or depth two less than its width.

Assume that there is a partition $\la$ in the decomposition of
$\mathrm{ind}_{S_{2k}}({\psi_{k}})$ that has a diagonal hook that does
not have the depth equal to one less than its width.  Assume that
$h=h_i$ is the smallest diagonal hook in $\la$ with $W(h) \neq D(h)+1$
and that $W(h) > D(h)+1$ (the case for $W(h) < D(h)+1$ follows
similarly).

First assume $h$ is the smallest diagonal hook in $\la$.  Let $h'$
be the hook constructed by removing the node $(i+D(h)-1,i)$ from
$h$.  Define $\la '$ to be the partition of $2k-1$ constructed by
replacing the hook $h$ in $\la$ by $h'$.  By the
Littlewood-Richardson Rule, $\la'$ is in the decomposition of
$\mathrm{res}_{2k-1}(\mathrm{ind}_{S_{2k}}({\psi_{k}}))$. Since the hook
$h'$ has $W(h) > D(h)+2$, it is not in the decomposition of
$\mathrm{ind}_{S_{2k-1}} (\mathrm{ind}_{S_{2k-2}}({\psi_{k-1}}))$.
Therefore, by Equation~\ref{eqn:equivdecomps}, $\la$ is not in the
decomposition of $\mathrm{ind}_{S_{2k}}({\psi_{k}})$.

If $h$ is not the smallest diagonal hook in $\la$, let $\la'$ be the
partition of $2k-1$ constructed by removing one cell from the smallest
diagonal hook in $\la$. This partition will not be in the decomposition of
$\mathrm{ind}_{S_{2k-1}} ( \mathrm{ind}_{S_{2k-2}}({\psi_{k-1}}))$
since it will have two diagonal hooks with width not one more than the depth.
Again by Equation~\ref{eqn:equivdecomps}, $\la$ is not in the
decomposition of $\mathrm{ind}_{S_{2k}}({\psi_{k}})$.
\qed

\begin{theorem}\label{thm:S2wrAkvalues}
The group $S_2\wr A_k $ is multiplicity free if and only if $k \in (3,4,7,8,11,12,16,20,24)$.
\end{theorem}
\proof
Assume $k=4a+1$ then all diagonal hooks in the following partition
have depth one more than its width $[(2a+2)^2, 2^{2a-1}]$. Since each
entry in the partition is even, this partition is in the decomposition
of both $\indg{S_2 \wr S_k}{S_{2k}}$ and $\mathrm{ind}_{S_2k}(\psi)$
so it occurs twice in the decomposition of $\indg{S_2 \wr
A_k}{S_{2k}}$.

Similarly, if $k=4a+2$ then the partition $[(2a+2)^2, 4, 2^{2a-2}]$
occurs twice in the decomposition of
$\indg{S_2 \wr A_k}{S_{2k}}$.

If $k=4a+3$ and $a\geq 4$ the partition $[2a,2a,6,6,6,2^{2a-3}]$
occurs twice in the decomposition of $\indg{S_2 \wr A_k}{S_{2k}}$. If
$a =3$ the partition $[6^5]$ occurs twice.  For $a \in [0,1,2]$ it is
possible to determine that $S_2 \wr A_k$ is multiplicity free either
by GAP or by checking the partitions of $2k$.

Finally, if $k=4a$ and $a \geq 7$ the partition
$[2a-6^2,8^5,2^{2a-14}]$ occurs twice in the decomposition of
$\indg{S_2 \wr A_k}{S_{2k}}$. For $a \leq 6$, by inspection of the
partitions in the decomposition of $\indg{S_2 \wr A_k}{S_{2k}}$ the
group $S_2 \wr A_k$ is multiplicity free.
\qed 

\begin{theorem}\label{thm:S2capA2kwrSk}
The group $((S_2)^k \cap A_{2k}) \rtimes S_k$ is multiplicity free if
and only if $k\in\{2,4,5,6,8,9,12,13,16,17,20,24,28,32\}$.
\end{theorem}
\proof
By Proposition~\ref{prop:linearcharsofwreath}, the decomposition of
$\indg{((S_2)^k \cap A_{2k}) \rtimes S_k }{S_2
\wr S_k}$ is $1+\sigma\psi_k$.  The decomposition of
$\mathrm{ind}_{S_{2k}}(\sigma \psi_k)$ contains the partitions which
are conjugate to the partitions in the decomposition of
$\mathrm{ind}_{S_{2k}}(\psi_k)$. {}From Theorem~\ref{thm:S2wrAk},
these are exactly the partitions in which all the diagonal hooks have width one
less than the depth.  For each value of $k$ not listed above, we will
give a partition of $2k$ in which all the diagonal hooks have width one less
than the depth and all the entries are even. Such a partition will
occur twice in the decomposition of $\indg{((S_2)^k \cap A_{2k})
\rtimes S_k }{S_{2k}}$ since it occur both in the decomposition of
$\mathrm{ind}_{S_{2k}}(1_{S_2 \wr S_k})$ and in the decomposition of
$\mathrm{ind}_{S_{2k}}(\psi_k)$.

For $k=4a$ and $a \geq 9$ the partition is $[(2a-10)^2,8^7,2^{2a-18}]$;
for $k=4a+1$ with $a \geq 5$ the partition is $[(2a-4)^2,6^5,2^{2a-10}]$;
for $k=4a+2$ with $a\geq 2$ the partition is $[(2a)^{2}, 4^3, 2^{2a-4}]$ and for $k=6$ the partition is $[2^3]$;
for $k=4a+3$ and all values of $a$ the partition is $[(2a+2)^{2}, 2^{2a+1}]$.

For the remaining cases, it is not hard to check that no partition of
$2k$ can satisfy the two requirements.
\qed

\begin{theorem} 
The group $((S_2)^k \cap A_{2k}) \rtimes A_k$ is not multiplicity free.
\end{theorem}
\proof
If $((S_2)^k \cap A_{2k}) \rtimes A_k$ is multiplicity free then each
of $(S_2 \wr S_k)\cap A_{2k}$, $S_2 \wr A_k$ and $(S_2)^k \cap
A_{2k}) \rtimes S_k$ is multiplicity free. There are no values of $k$ for
which all three of these groups are multiplicity free.
\qed

\subsubsection{Multiplicity free subgroups of $S_\ell \wr S_3$ and $S_3 \wr S_k$}

In this section we consider the remaining four structures for the
blocks of imprimitivity for imprimitive multiplicity-free
groups. These structures are three blocks of size three, four or five
or four blocks of size three. 

Let $\ell$ be the size of the blocks of imprimitivity.  By
Corollary~\ref{cor:thomonblocks}, any imprimitive multiplicity-free
subgroup $G$ of $S_\ell \wr S_3$ must be $t$-homogeneous on the blocks
of imprimitivity where $t \leq \ell$. If $\ell \neq 5$, then the group
$G$ must contain the group $A_\ell^3$.  Using the same method as in
Theorem~\ref{thm:subgroupsSkwrS2}, the group $G/A_\ell^3$ must be a
multiplicity-free subgroup of $(S_\ell \wr S_3)/ A_\ell^3 \cong S_2
\wr S_3$.  {}From the previous section, we know all the
multiplicity-free subgroups of $S_2 \wr S_3$.  First we find the
groups $G$, which have the property that $G/A_\ell^3$ is one of the
multiplicity-free subgroups of $S_2 \wr S_3$. Next, using GAP, we
determine which of these groups are multiplicity free in
$S_{3\ell}$. 

Define a subgroup $D_\ell$ of $S_\ell^3$ by 
\[
D_\ell=\{ (a,b,c) \in S_\ell^3 : \sgn(a) =\sgn(b)=\sgn(c)\}.
\]
The group $D_\ell$ is normal in $S_\ell \wr S_3$, so $SD_\ell=D_\ell \rtimes
S_3$ is a subgroup of $S_\ell \wr S_3$. Further, define a subgroup of $SD_\ell$ by
\[
RD_{\ell} =\{ (a,b,c:\pi) \in S_\ell \wr S_3 : \sgn(\pi) =\sgn(a)=\sgn(b)=\sgn(c)\}.
\]
If $\ell$ is odd then $RD_{\ell} = SD_\ell \cap A_{3\ell}$. 

The following table lists all multiplicity-free subgroups of $S_\ell
\wr S_3$, the columns are headed by the ``supergroup''.

\renewcommand{\arraystretch}{1.25} 
\begin{center}
\begin{tabular}{l|l|l}
$S_3 \wr S_3$ & $S_3 \wr S_4$ & $S_4 \wr S_3$\\  \hline
$(S_3 \wr S_3) \cap A_9$ & $(S_3 \wr S_4) \cap A_{12}$ & $(S_4 \wr S_3) \cap A_{12}$ \\
$(S_3^3 \cap A_9) \rtimes S_3$ & &$(S_4^3 \cap A_{12}) \rtimes S_3$\\
 & & $SD_4$ \\
 & & $RD_4$ \\
\end{tabular}
\end{center}
\renewcommand{\arraystretch}{1}

A similar method can be used to determine the multiplicity-free
subgroups of $S_5 \wr S_3$.
\begin{theorem}
The transitive, subgroups of $S_5 \wr S_3$ 
which are multiplicity free in $S_{15}$ are
\begin{enumerate}[(a)]
\item $(S_5 \wr S_3)\cap A_{15}$,
\item $A_5 \wr S_3$,
\item $((S_5)^3 \cap A_{15}) \rtimes S_3$,
\item $SD_5$,
\end{enumerate}
\end{theorem}
\proof
If $G$ is a multiplicity-free subgroup of $S_5 \wr S_3$ then either
$A_\ell^3\leq G$ or $AGL(1,5)^3 \leq G$.

For the first case, we only need to check which of the groups $G$ with
$G/A_\ell^3$ multiplicity free in $S_2 \wr S_3$ are multiplicity free
in $S_{15}$.  With a calculation in GAP, the groups $(S_5 \wr S_3)\cap
A_{15}$ and $A_5 \wr S_3$ are multiplicity free. Since $A_5 \wr
S_3$ is a subgroup of $((S_5)^3 \cap A_{15}) \rtimes S_3$ and $SD_5$,
these two groups are also multiplicity free.

In the second case, $G$ must be a subgroup of $AGL(1,5) \wr S_3$. Since
this group is not multiplicity free, $G$ is not multiplicity free.
\qed

\section{Intransitive Subgroups}

Let $\Omega$ be a set of cardinality $n$.  Throughout this section we
assume that the group $G \leq S_n$ is intransitive on $\Omega$. For a
group $G \leq S_n$, we will use $S_1\times G$ to denote the subgroup
of $S_{n+1}$ that fixes one point of the underlying set and is $G$ on
the remaining $n$ points.

{}From Lemma~\ref{lem:saxl}, if $G$ is multiplicity free and transitive, then it has
exactly two orbits on $\Omega$. Denote these two orbits by $\Gamma$ and
$\Delta$. Then $G \leq G_{\upharpoonright\Gamma} \times
G_{\upharpoonright\Delta}$.

\begin{lemma}\label{lem:orbitshom}
Let $G$ be an intransitive, multiplicity-free subgroup of $S_n$.  Let
$\Gamma$ and $\Delta$ denote the two orbits of $G$ and assume
$|\Gamma|=k$ and $|\Delta| =n-k$ with $2k \leq n$.  Then
\begin{enumerate}[(a)]
\item{for all $\ell$, with $1\leq \ell \leq k$, the group $G$ has exactly $\ell+1$ orbits on
           ${\Omega \choose \ell}$;}
\item{for all $\ell \leq k$, the groups $G_{\upharpoonright\Gamma}$ and
      $G_{\upharpoonright\Delta}$ are $\ell$-homogeneous.}
\end{enumerate}
\end{lemma}

\proof
First we prove part (a). Let $\ell$ be an integer with $1 \leq \ell
\leq k$.  Since $G$ is multiplicity free, by Lemma~\ref{lem:saxl},
there can be at most $\ell+1$ orbits on ${\Omega \choose \ell}$.
To see that there are exactly this many orbits on ${\Omega \choose
\ell}$, consider the intersection of any $\ell$-set from $\Omega$ with
the set $\Gamma$.  The size of this intersection can take $\ell+1$
different values, 0 to $\ell$. Since $\Gamma$ is an orbit for $G$, if
the size of the intersection with $\Gamma$ for two $\ell$-sets from
$\Omega$ is different, then the two $\ell$-sets must be in different
orbits of $G$.  Hence $G$ exactly $\ell+1$ orbits on ${\Omega \choose
\ell}$.

Next we prove part (b). The proof of part (a) implies that for each
integer $\ell$ with $0\leq \ell \leq k$, there is exactly one orbit of
$G$ on ${\Omega \choose k}$ where the intersection of the $k$-sets in
the orbit with $\Gamma$ has cardinality $\ell$.  In particular, for
any $\ell \leq k$, the restricted group $G_{\upharpoonright\Gamma}$
has a single orbit on ${\Gamma \choose \ell}$ and
$G_{\upharpoonright\Gamma}$ is $\ell$-homogeneous for every $\ell \leq
k$.

Similarly $G_{\upharpoonright\Delta}$ is $\ell$-homogeneous for every
$\ell \leq k$.
\qed

\begin{lemma}\label{lem:productsmf}
Let $G$ be an intransitive, multiplicity-free subgroup of $S_n$.
Assume that the two orbits of $G$ on $\Omega$ are $\Gamma$ and $\Delta$ with
cardinality $k$ and $n-k$ respectively.  Then
$G_{\upharpoonright\Gamma}$ and $G_{\upharpoonright\Delta}$ are
multiplicity-free subgroups of $S_k$ and $S_{n-k}$.
\end{lemma}
\proof
Assume that
\[
\indg{G_{\upharpoonright\Gamma}}{S_k} = \sum_{\mu \vdash k}a_\mu V_\mu
\]
and 
\[
\indg{G_{\upharpoonright\Delta}}{S_{n-k}} = \sum_{\nu \vdash n-k}b_\nu V_\nu.
\]
By the Littlewood-Richardson Rule,
\[
\indg{G_{\upharpoonright\Gamma} \times G_{\upharpoonright\Delta}}{S_n}
	= \sum_{\mu \vdash k}\sum_{\nu \vdash n-k}
	\sum_{\la \vdash n}a_\mu b_\nu N_{\mu,\nu}^\la V_\la.  
\]
If either $a_\mu> 1$ for some $\mu \vdash k$ or $b_\nu>1$ for some
$\nu \vdash (n-k)$ then $G_{\upharpoonright\Gamma} \times
G_{\upharpoonright\Delta}$ is not multiplicity free. Since $G \leq
G_{\upharpoonright\Gamma} \times G_{\upharpoonright\Delta}$, by
Theorem~\ref{thm:basicfacts}(\ref{subgroups}), $G$ would also
not be multiplicity free.
\qed

Lemma~\ref{lem:orbitshom} and Lemma~\ref{lem:productsmf} impose strong
restrictions on the groups $G_{\upharpoonright\Gamma}$ and
$G_{\upharpoonright\Delta}$. In
particular, the group $G_{\upharpoonright\Gamma}$ must be either the
symmetric group, the alternating group, $AGL(1,5)$ with $k=5$,
$PGL(2,5)$ with $k=6$ or $P\Gamma L(2,8)$ with $k=9$.  There are more choices for
$G_{\upharpoonright\Delta}$, since it is only $\ell$-homogeneous for
$\ell \leq |\Gamma| \leq |\Delta|$.  But as
$G_{\upharpoonright\Delta}$ must be transitive and multiplicity free,
$G_{\upharpoonright\Delta}$ must be one of the groups from
Section~\ref{sec:transitive}.

We will first determine which products of groups from
Section~\ref{sec:transitive} are multiplicity free and then consider
the subgroups of these product groups.  Proposition~\ref{prop:StimesS}
lists the multiplicity-free groups which are the product of only
symmetric and alternating groups.

\begin{proposition}\label{prop:specialproducts}
Let $k$ be an integer, then
\begin{enumerate}[(a)]
\item $AGL(1,5) \times S_k$, $PGL(2,5) \times S_k$ and $P\Gamma L(2,8) \times S_k$ are multiplicity free 
             for all positive integers $k$,
\item  $AGL(1,5) \times A_k$ is multiplicity free if and only if $k \geq 4$,
\item  $PGL(2,5) \times A_k$ is multiplicity free if and only if $k \geq 3$.
\item  $P\Gamma L(2,8) \times A_k$ is multiplicity free for any $k=7$ or $k >=11$.
\end{enumerate}
\end{proposition}

\proof
By the Littlewood-Richardson rule, the decomposition into irreducible
characters of $AGL(1,5) \times S_1$ considered as a subgroup of $S_6$ is
\[
[6]+[5,1]+[3,2,1]+[2^3]+ [2^2,1^2]
\]
and for all integers $k> 1$, the decomposition into irreducible
characters of $AGL(1,5) \times S_k$ is
\[
\sum_{i=0}^{\min(k,\lfloor\frac{k+5}{2} \rfloor)} [k+5-i,i]
       + [k+2,2,1]+[k+1,2^2]+[k+1,2,1^2]+[k,2^2,1].
\]

Similarly, the decomposition into irreducible
characters of $PGL(2,5) \times S_1$ as a subgroup of $S_7$ is
\[
[7]+[6,1]+[3,2^2]+[2^3,1]
\]
and the decomposition of $PGL(2,5) \times S_k$ is
\[
\sum_{i=0}^{\min(k,\lfloor\frac{k+6}{2} \rfloor)}[k+6-i,i]
       +[k+2,2^2]+[k+1,2^2,1]+[k,2^3].
\]

For $k \leq 4$ it can be checked using GAP that $P\Gamma L(2,8) \times S_k$ is multiplicity free.
For $k > 4$ the decomposition into irreducible
characters of $P\Gamma L(2,8) \times S_k$ is
\begin{align*}
&\sum_{i=0}^{\min(k,\lfloor\frac{k+9}{2} \rfloor)} [k+9-i,i] + [k+1,1^8]+[k,1^9] \\
& + \sum_{i=0}^{4} \left( [5+k-i,i+1,1^3] + [4+k-i,i+1,1^4]\right. \\
& \left. +[4+k-i,4,i+1] + [3+k-i,4,i+1,1] \right) \\
& + [3+k,2,2,2] +[k+3,2,2,2,2]+[k+2,3,2,2,2] \\
& + [3+k-1,3,2,2] + [3+k-1,2,2,2,1] + [3+k-2,3,2,2,1] 
\end{align*}

The decomposition into irreducible characters of $AGL(1,5) \times A_k$ is
\begin{align*}
  &\sum_{i=0}^{\min(k,\lfloor\frac{k+5}{2} \rfloor)} [k+5-i,i]
   + [k+2,2,1]+[k+1,2^2]+[k+1,2,1^2] \\
  &+ [k,2^2,1]+[6,1^{k-1}] +[5,1^k] + [3^2,2,1^{k-3}] +[3^2,1^{k-1}] \\
  &+[3,2^2,1^{k-2}] +[3,2,1^{k}] +[2^3,1^{k-1}]+[2^2,1^{k+1}].
\end{align*}
If $k \geq 4$, then these irreducible characters are all different.

The decomposition into irreducible characters of $PGL(2,5) \times A_k$ is
\begin{align*}
&\sum_{i=0}^{\min(k,\lfloor\frac{k+6}{2} \rfloor)} [k+6-i,i]
       + [k+2,2^2]+[k+1,2^2,1]+[k,2^3]+[7,1^{k-1}] \\
&+ [6,1^k]+[3^3,1^{k-3}]+[3^2,2,1^{k-2}]+[3,2^2,1^{k-1}] +[2^3,1^{k}].
\end{align*}
If $k \geq 3$ all these irreducible characters are distinct.

For $k<11$ the decomposition into irreducible characters for
$P\Gamma L(2,8) \times A_k$ can be calculated in GAP.  If $k\geq 11$ Then the
the decomposition of $P\Gamma L(2,8) \times A_k$ is the sum of the
partitions formed by adding $k$ boxes to the partitions in the
decomposition of of $P\Gamma L(2,8)$ first, in such a way that no two boxes
are in the same row and second in such a way that no two boxes are in the same column.
The decomposition of $P\Gamma L(2,8)$ is 
\[
[9] + [5,1^4] + [4,4,1] + [3,2,2,2]  + [1^9]
\]
it is tedious but not difficult to see that for $k\geq 11$ each of
these partitions is distinct.
  \qed

\begin{proposition}\label{prop:productwithwreath}
Let $k$, $2\leq c$ and $2 \leq d$ be integers and let $G$ be a subgroup of $S_k$.
Then $G \times (S_c \wr S_d)$ is multiplicity free if and only if
$G=S_1$ and $c=2$ or $G=S_1$ and $d=2$.
\end{proposition}
\proof
If $G \times (S_c \wr S_d)$ is multiplicity free, then it has exactly
two orbits on $\{1, \dots, cd+k\}$.  The first orbit is $\Gamma =\{1,
\dots ,k\}$ and the second is $\Delta=\{k+1,\dots, cd+k\}$.  By
Lemma~\ref{lem:orbitshom}(b), if $G \times (S_c \wr S_d)$ is
multiplicity free, then $S_c \wr S_d$ is $k$-homogeneous. Since for
all integers $c\geq 2$ and $d\geq 2$, the group $S_c \wr S_d$ is
imprimitive it is not $k$-homogeneous for any $k>1$.  Thus $k=1$ and
$G$ is $S_1$.

If $S_1 \times (S_c \wr S_d)$ is multiplicity free, then $S_c \wr S_d$
is also multiplicity free.  The only integers $3 \leq c$ and $3\leq d$
with $S_c \wr S_d$ multiplicity free are $(3,3)$, $(3,4)$, $(3,5)$ and
$(4,3)$. For each of these values of $(c,d)$, the decomposition of
$\indg{S_c \wr S_d}{S_{cd}}$ into irreducible characters includes the characters
$[cd-2,2]$ and $[cd-3,3]$.  By the Littlewood-Richardson
rule, $[cd-2,2] \circ [1] = [cd-1,2]+[cd-2,3]$ and 
$[cd-3,3] \circ [1] = [cd-2,3] + [cd-3,2]$.  For $c \geq 2$ and $d\geq
2$, the irreducible representation $[cd-2,3]$ occurs with multiplicity
at least two in the decomposition of $\indg{S_1 \times (S_c \wr
S_d)}{S_{cd+1}}$. Thus $c=2$ or $d=2$.

For a partition $(\la_1,\la_2, \dots, \la_a)$ of $d$ let
$2\la+1$ denote all the partitions of $2d+1$ of the form
$(2\la_1,2\la_2, ... ,2\la_i+1,... ,2\la_a)$, where
$\la_{i-1} > \la_{i}$.
Then Proposition~\ref{prop:wreath1}, together with the Littlewood-Richardson Rule, gives
that
\[
\indg{S_1 \times (S_2 \wr S_d)}{S_{2d+1}} = \sum_{\la \vdash c} 2\la+1.
\] 
Each of these characters is distinct and $S_1 \times (S_2 \wr S_d)$ is
multiplicity free. The number of characters in the decomposition is
the sum over all partitions $\la =(\la_1, \la_2,
\dots, \la_a)$ of $d$ of the number of $i \in \{0,\dots,a\}$ with $\la_{i} >
\la_{i+1}$ where we set $\la_0 =d+1$ and $\la_{a+1}=0$.

The irreducible decomposition of $\indg{S_1 \times (S_c \wr S_2)}{S_{2c+1}}$ is
\[
  \sum_{i=0}^{\lfloor c/2\rfloor}
     \left ( [2c-2i+1,2i]+[2c-2i,2i+1]+[2c-2i,2i,1] \right ).
\]
Since these characters are distinct, $S_1 \times (S_c \wr S_2)$ is multiplicity free.
\qed

For each multiplicity-free subgroup of $S_\ell \wr S_2$ the
decomposition of the induced representation is given in the proof of
either Theorem~\ref{thm:subgroupsSkwrS2} or the proof of
Theorem~\ref{thm:specialwreaths}. So for each subgroup, using the
Littlewood-Richardson rule, it is possible to determine if the direct
product of $S_1$ and the subgroup is multiplicity free.

\begin{proposition}
\begin{enumerate}[(a)]
\item The group $S_1 \times ( ( S_\ell\wr S_2 )\cap A_{2\ell})$ is multiplicity free for $\ell \geq 3$.
\item The group $S_1 \times (((S_\ell \times S_\ell) \cap A_{2\ell})\rtimes S_2)$ is multiplicity free for $\ell \geq 3$.
\item For all values of $\ell$, the group
$S_1 \times (A_\ell \wr S_2)$ is not multiplicity free.
\item The group $S_1 \times (AGL(1,5) \wr S_2)$ is not multiplicity free.
\item The group $S_1 \times (PGL(2,5) \wr S_2)$ is not multiplicity free. \qed
\end{enumerate}
\end{proposition}

\begin{proposition}
The groups $S_1 \times ((S_2 \wr S_k) \cap A_{2k})$, 
$S_1 \times (S_2 \wr A_k)$, and 
$S_1 \times (((S_2)^k \cap A_{2k}) \rtimes S_k)$ are not multiplicity free.
\end{proposition}
\proof
For all $k$, the partitions $[2^k]$ and $[2^{k-1},1^2]$ both occur in
the decomposition of $\indg{(S_2 \wr S_k) \cap A_{2k}}{S_{2k}}$. by the
Littlewood-Richardson rule the partition $[2^{k},1]$ occurs twice in the decomposition of
$\indg{S_1 \times ((S_2 \wr S_k) \cap A_{2k})}{S_{2k+1}}$

The fact that the groups $S_1 \times (S_2 \wr A_k)$, and 
$S_1 \times (((S_2)^k \cap A_{2k}) \rtimes S_k)$ are not multiplicity free can be
tested using GAP for the list of values of $k$ given in
Theorems~\ref{thm:S2wrAkvalues} and \ref{thm:S2capA2kwrSk}.
\qed

All other product groups are given in Table~\ref{table:intrans}. These
group were determined to be multiplicity free using GAP.

Finally we consider which proper subgroups of the multiplicity-free
product groups are also multiplicity free. A subgroup $G$ of a product
group $H \times K$ is a \textsl{subdirect product} if for $\pi_1$ and
$\pi_2$ the projection maps on to the first and second factor,
$\pi_1(G) = H$ and $\pi_2(G)=K$. We give a well-known characterization
of subdirect products. This result is standard but we include that
proof for completeness.

\begin{proposition}
Let $G$ be a subdirect product of $H \times K$.  Then for some
pair of homomorphisms $\phi$ of $H$ and $\psi$ of $K$ both to the same group 
\[
G=\{(h,k) \in H \times K : \phi(h)=\psi(k)\}.
\]
\end{proposition}

\proof
Define 
\[
H_1 =\{h : (h,1) \in G\}, \quad K_1 = \{k : (1,k) \in G\}.
\]
The subgroup $H_1 \times 1$ is the kernel of the projection of $G$ to
$K$, thus $H_1 \times 1$ is normal in $G$.  The image of the projection of $H_1
\times 1$ in $H$ is $H_1$ and since the projection map is a homomorphism,
$H_1$ is a normal subgroup of $H$.  Similarly, $K_1$ is a normal
subgroup of $K$.

Since $1 \times K_1$ is the kernel of the surjective projection of $G$
onto $H$, the group $G /(1 \times K_1)$ is isomorphic to
$H$. Similarly $G /(H_1 \times 1) \cong K$.

Then 
\[
\frac{H}{H_1} \cong
\frac{G/(1 \times K_1)}{H_1 \times 1} \cong 
\frac{G}{(H_1 \times 1)(1 \times K_1)}\cong 
\frac{G/(H_1 \times 1)}{1 \times K_1} \cong
\frac{K}{K_1}. 
\] 
Thus $H/H_1 \cong K/K_1$.

Let $\rho:H/H_1 \rightarrow K/K_1$ be an isomorphism. Define $\phi'$ be
the natural homomorphism of $H$ to $H/H_1$ and $\phi = \rho
\circ \phi'$. Let $\psi'$ be the natural
homomorphism of $K$ to $K/K_1$. We will show that $G = \{(h,k) \in H
\times K : \phi(h)=\psi(k)\}$.  

Define the map $\Psi:H \times K \rightarrow K/K_1$ by $\Psi(h,k) =
\phi(h)\psi(k^{-1})$. The claim is that $G$ is the kernel of $\Psi$.
If $(h,k) \in G$ then $\Psi(h,k) =\phi(h)\psi(k^{-1}) = 1$. Thus $G \subseteq \ker(\Psi)$.
Since
\[
\frac{|H \times K|}{ |\ker(\Psi)|} = |K/K_1| = \frac{|H||K|}{|G|}.
\]
Thus $|G| = |\ker(\Psi)|$ and $G = \ker(\Psi)$.
\qed

For all of the intransitive multiplicity-free, product groups, at
least one of the products is either the symmetric group or the
alternating group. With the exception of $k=4$, the only normal
subgroup of $S_k$ is the alternating group and $A_k$ has no non-trivial normal
subgroups.  The only non-trivial homomorphism of $S_k$ is the sign map. For $k\neq
4$ and a group $H$, if $G$ is a proper subgroup of $S_k \times H$ then
\[
G=\{(h,k): \sgn(h)=\sgn(k)\} =(S_k \times H) \cap A_n.
\]
For each multiplicity-free intransitive group, we will
determine if its intersection with the alternating group is also multiplicity free.
There are no non-trivial homomorphisms of $A_k$, so we do not need to check
subgroups of $A_k \times H$.  Finally we will consider the special
cases of proper subgroups of $S_4 \times H$ and $A_4 \times H$.

\begin{proposition}\label{prop:StimesScapA}
For all positive integers $n$ and $k$, except $k=2$ and $n=4$, the group $(S_k \times
S_{n-k})\cap A_{n}$ is multiplicity free.
\end{proposition}
\proof
The decomposition of $(S_k \times S_{n-k})\cap A_{n}$ into irreducible
characters is 
\[
\sum_{i=0}^{ \min\{k,n-k\} } \left([n-i,i] + [2^i,1^{n-2i}]\right).\qed
\]

\begin{proposition}
For all integers $k \geq 4$, the group $(AGL(1,5) \times S_k)\cap
A_{k+5}$ is multiplicity free.  For all positive integers $k$ the
group $(PGL(2,5) \times S_k)\cap A_{k+6}$ is multiplicity free.
\end{proposition}
\proof
By Proposition~\ref{prop:intersectalt} and \ref{prop:specialproducts},
the decomposition into irreducible characters of $(AGL(1,5) \times
S_k) \cap A_{k+5}$ and $k>1$ is
\begin{align*}
& \sum_{i=0}^{\min(k,\lfloor\frac{k+5}{2} \rfloor)} \left( [k+5-i,i]+[2^i,1^{k+5-2i}]\right)
       + [k+2,2,1] \\
     &+ [k+1,2^2]    +[k+1,2,1^2]  +[k,2^2,1] + [3,2,1^k] \\
     &+ [3^2,1^{k-1}]+[4,2,1^{k-1}]+[4,3,1^{k-2}].
\end{align*}
For $k\geq 4$ each of these characters is unique. Further, $(AGL(1,5)
\times S_1) \cap A_{6}$ is not multiplicity free.

The decomposition into irreducible
characters of $(PGL(2,5) \times S_k) \cap A_{k+6}$ for $k=1$ is
\[
[7]+[1^7]+[6,1]+[2,1^5]+[3,2^2]+[3^2,1]+[2^3,1]+[4,3]
\]
and in general it is
\begin{align*}
& \sum_{i=0}^{\min(k,\lfloor\frac{k+6}{2} \rfloor)} \left( [k+6-i,i]+ [2^i,1^{k+6-2i}] \right)
       + [k+2,2^2]+[k+1,2^2,1] \\
      & +[k,2^3] + [3^2,1^{k}]+[4,3,1^{k-1}]+[4^2,1^{k-2}]. 
\end{align*}
For all positive integers $k$ each of these irreducible representations are distinct.
\qed

For no other groups from Table~\ref{table:intrans} is the
intersection with the alternating group is multiplicity free.  This
can be tested using GAP or by checking the decomposition for conjugate
partitions.

Finally, we need to determine which proper subgroups of the
multiplicity-free groups of the form $S_4 \times G$ and $A_4 \times G$
are multiplicity free.  These subgroups are special since the groups
$S_4$ and $A_4$ have homomorphisms both to $\mathbb{Z}_2$ and $S_3$.

The subgroups of the form $(S_4 \times G) \cap A_n$ and $(A_4 \times
G) \cap A_n$ have been determined, so we only need to consider groups
$G$ that have a homomorphism to $S_3$. If $G$ has such a homomorphism
then it must have a normal group with index dividing 6.  {}From the list of
multiplicity-free product groups we only need to consider the
following: $S_4 \times S_4$ and $S_4 \times A_4$.  Using GAP to check
all subgroups of these groups we find that the only proper
multiplicity-free subgroup is $(S_4 \times S_4) \cap A_8$.

\section{Conclusion and Further Work}

In our introduction we defined a set of matrices which are the
adjacency matrices of the graphs whose vertex set is the collection of
all uniform $k$-partitions of an $n$-set with adjacencies defined by
the meet tables. We know precisely for which values of $k$ and $n$
this set of matrices commute. The adjacency matrix of graph $QI(12,3)$
is one of these matrices. The eigenvalues of this matrix are
known~\cite{karensthesis} and using standard eigenvalue methods we get
that the size of the maximum clique of $QI(12,3)$ is 7.

The next stage of this project is to develop tools for the
non-commutative case; we want better bounds on the size of the maximum
cliques in $QI(n,k)$ and we want to know if these graphs are cores.
Since $S_k \wr S_k$ is not multiplicity free when $k >3$, the method
of Godsil and Newman~\cite{Godsil:Newman} will not extend directly.

In Proposition~\ref{prop:linearcharsofwreath} we showed that the group
$S_\ell \wr S_2$ has exactly three linear characters, not including
the trivial character.  It is interesting that the corresponding
induced characters of $S_{2\ell}$ are multiplicity free.  This
suggests the problem of determining which subgroups of the symmetric
group have a linear character whose induced character is multiplicity
free.  (It is not clear that solving this would have any combinatorial
impact, but nonetheless it could be interesting.)

Finally, for each of the multiplicity-free subgroups we found,
there is a corresponding association scheme with the cosets of the
subgroup as its vertices.  Some of these may prove interesting in
their own right.  We plan to search for fusion schemes of low rank in
these schemes.

\section{Acknowledgements}
We would like to thank Mark Wildon who pointed out an error in Proposition~\ref{prop:specialproducts}.

\newpage
\section{Tables of Multiplicity-free Groups}
For three families of groups ( $S_2 \wr A_k$, $(S_2^k \cap
A_{2k})\rtimes S_k$ and $S_1 \times (S_2 \wr S_k)$ ) the rank is too
complicated to be included in the table.

\begin{table}[h]
\begin{center}
\begin{tabular}{l ccc}
Group & $n$ & index & rank\\
$A(n)$ & $n$ & $2$ & 2 \\
$S(n)$ & $n$ & $1$ & 1 \\
$AGL(1,5) \cap A_5$ & 5 & 12 & 4\\
$AGL(1,5)$ & 5 & 6 & 2 \\
$PSL(2,5)$ & 6 & 12 & 4 \\
$PGL(2,5)$ & 6 & 6 & 2 \\
$AGL(1,7)$ & 7 & 120 & 7 \\
$PSL(3,2)$ & 7 & 30 & 4 \\
$A\Gamma L(1,8)$ & 8 & 240 & 8 \\
$PGL(2,7)$ & 8 & 120 & 5 \\
$AGL(3,2)$ & 8 & 30 & 4 \\
$AGL(2,3)$ & 9 & 840 & 9 \\
$P\Gamma L(2,8)$ & 9 & 240 & 5 \\
$\aut(PSL(2,9))$ & 10 & 2520 & 10 \\
$M_{11}$ & 11 & 5040 & 10 \\
$M_{11}$ & 12 & 60480 & 26 \\
$M_{12}$ & 12 & 5040 & 8 \\
\end{tabular}
\caption{Multiplicity-free transitive primitive groups. \label{table:transitive}}
\end{center}
\end{table}

\renewcommand{\arraystretch}{1.25} 
\begin{table}
\begin{center}
\begin{tabular}{lccc}
Group & $n$ & index & rank\\
$S_\ell \wr S_2$ & $2\ell$ & $\frac{1}{2}{2\ell \choose \ell}$  & ${\lfloor \ell/2 \rfloor }+1$ \\
$(S_\ell\wr S_2)\cap A_{2\ell}$, $\ell \geq 3$ & $2\ell$ & ${2\ell \choose \ell}$ &  $2(\lfloor \ell/2 \rfloor +1 )$ \\
$((S_\ell \times S_\ell) \cap A_{2\ell})\rtimes S_2$, $\ell \geq 3$ &$2\ell$ &${2\ell \choose \ell}$ & $2(\lfloor \ell/2 \rfloor +1 )$ \\
$A_\ell \wr S_2$, $\ell \geq 3$  & $2\ell$ & $2{2\ell \choose \ell}$ & $2(\lfloor \ell/2 \rfloor +1 ) +2$\\
$AGL(1,5) \wr S_2$ & 10 & 4536 & 16 \\
$PGL(2,5) \wr S_2$ & 12 &  16632 & 13 \\

$S_2\wr S_k$  & $2k$ & $\frac{(2k)!}{2^k k!}$ & $p(k)$ \\
$(S_2\wr S_k) \cap A_{2k}$, $k$ odd & $2k$ & $\frac{(2k)!}{2^{k-1} k!}$ & $2p(k)$ \\ 

$S_2 \wr A_k$ & $2k$ & $\frac{(2k)!}{2^{k-1} k!}$ &  \\
\ ($k \in (3,4,7,8,11,12,16,20,24)$) & & & \\

$(S_2^k \cap A_{2k}) \rtimes S_k$ & $2k$ & $\frac{(2k)!}{2^{k-1} k!}$ &  \\
\ ($k\in\{2,\! 4,\!5,\!6,\!8,\!9,\!12,\!13,\!16,\!17,\!20,\!24,\!28,\!32\}$)& & & \\

$S_3 \wr S_3$ & 9 & 280 & 5 \\
$(S_3 \wr S_3) \cap A_9$ & 9 & 560  & 10 \\
$((S_3)^3 \cap A_9) \rtimes S_3$ & 9 & 560  & 9 \\
$S_3 \wr S_4$ & 12 & 15400  & 12 \\
$(S_3 \wr S_4) \cap A_{12}$ & 12 & 30800 & 24 \\
$S_4 \wr S_3$ & 12 & 5775 & 9 \\
$(S_4 \wr S_3) \cap A_{12}$& 12 & 11550 & 18\\
$((S_4)^3 \cap A_{12}) \rtimes S_3$& 12 & 11550 & 16 \\
$SD_4$ & 12 & 23100& 18 \\
$RD_4$ & 12 & 46200 & 33\\
$S_5 \wr S_3$ & 15 & 126126 & 13 \\
$(S_5 \wr S_3)\cap A_{15}$ & 15 &252252 & 26\\
$((S_5)^3 \cap A_{15}) \rtimes S_3$ & 15 & 252252 & 25\\
$SD_5$ & 15 &504504 & 24\\
$A_5 \wr S_3$ & 15 & 1009008 & 46 \\
\end{tabular}
\caption{Multiplicity-free transitive imprimitive groups\label{table:imprimitive}
}
\end{center}
\end{table}

\begin{table}
\begin{center}
\begin{tabular}{lccc}
Group & $n$ & index & rank\\
$S_1 \times AGL(1,5)$ & $6$ & 36 & 5 \\
$S_1 \times PSL(2,5)$  & 7 & 84 & 8\\
$S_1 \times PGL(2,5)$  & 7 & 42 & 4 \\
$S_1 \times PSL(3,2)$ &  8 & 240 & 10\\
$S_1 \times PGL(2,7)$ &  9 & 1080 & 12\\
$S_1 \times AGL(3,2)$ &  9 & 270 & 8 \\
$S_1 \times P\Gamma L(2,8)$  & 10 & 2400 & 13 \\
$S_1 \times M_{12}$ &  13& 65520 &18 \\
$S_2 \times PSL(2,5))$  & 8 & 336 & 11\\
$S_2 \times PSL(3,2)$  & 9 & 1080 & 14 \\
$S_2 \times AGL(3,2)$  & 10 & 1350 & 11\\
$S_2 \times P\Gamma L(2,8)$ & 11 & 13200 & 20\\
$S_2 \times M_{12}$ & 14 & 458640 & 28 \\
$S_3 \times AGL(3,2)$ & 11 & 4950 & 13\\
$A_3 \times AGL(3,2)$ & 11 & 9900 & 26 \\
$S_3 \times P\Gamma L(2,8)$ & 12 & 52800 & 26 \\
$S_3 \times M_{12}$ & 15 & 2293200 & 37 \\
$S_4 \times P\Gamma L(2,8)$ & 13& 171600 & 30 \\
$S_4 \times M_{12}$ & 16 & 9172800 & 44 \\
$S_5 \times M_{12}$ & 17 & 31187520 & 48 \\
\end{tabular}
\caption{Multiplicity-free intransitive groups\label{table:intrans}}
\end{center}
\end{table}

\begin{table}
\begin{center}
\begin{tabular}{llccc}
Group & & $n$ & index & rank\\
$S_k \times S_{n-k}$ & $(2k \leq n)$ & $n$ & ${n \choose k}$ & $k+1$ \\
$(S_k \times S_{n-k})\cap A_{n}$ & $(2k\leq n$, $(k,n)\neq (2,4))$& $n$ & $2{n \choose k}$ & $2k+2$ \\
$A_k \times S_{n-k}$ & $(2k \leq n$, $k \neq 2)$ & $n$ & $2{n \choose k}$ &  $k+3$\\
$S_k \times A_{n-k}$ & $(2k \leq n$, $k \neq n-2)$ & $n$ & $2{n \choose k}$ &  $k+3$\\
$A_k \times A_{n-k}$ & $(k\geq 3$, $2k \leq n-2)$ & $n$ & $4{n \choose k}$ &  $2k+6$\\

$S_k \times AGL(1,5)$ & $(k\geq 2)$ & $k+5$ & $\frac{(k+5)!}{20 (k!)}$ & $k+5$ \\
$A_k \times AGL(1,5)$ & $(k\geq 4)$ & $k+5$ & $\frac{(k+5)!}{10 (k!)}$ & $k+13$ \\
$(S_k \times (AGL(1,5)) \cap A_{k+5}$ &$(k\geq 5)$ & $k+5$ & $\frac{(k+5)!}{10 (k!)}$ & $2k+10$ \\

$S_k\times PGL(2,5) $ & $(k\geq 2)$ & $k+6$ & $\frac{(k+6)!}{120 (k!)}$ & $k+4$ \\
$A_k \times PGL(2,5) $ & $(k\geq 3)$ & $k+6$ & $\frac{(k+6)!}{60 (k!)}$ & $k+10$ \\
$(S_k \times (PGL(2,5)) \cap A_{k+6}$& $(k\geq 2)$ & $k+6$ & $\frac{(k+6)!}{60 (k!)}$ & $2k+8$ \\

$S_k\times P\Gamma L(2,8) $ & $(k\geq 5)$ & $k+9$ & $\frac{240 (k+9)!}{9!k!}$ 
& \small{$27\!+\!\min\{ k, \lfloor \frac{k+9}{2} \rfloor\}$}\\
$A_k\times P\Gamma L(2,8) $ & $k=7$ or $k>=11$ & $k+9$ & $\frac{480 (k+9)!}{9!k!}$ & 
\small{$52+k$ when $k\geq 16$}\\

$S_1 \times (S_k \wr S_2)$ & &$2k+1$ & $\frac{(2k+1)!}{2(k!)^2}$ & $3\lfloor \frac{k}{2}\rfloor$\\
$(S_1\times((S_k \wr S_2))\cap A_{2k+1})$& $(k \geq 3)$&$2k+1$ &$\frac{(2k+1)!}{(k!)^2}$&$6\lfloor \frac{k}{2}\rfloor$\\
$(S_1\times ((S_k^2 \cap A_{2k})\rtimes S_2)$ & $(k \geq 3)$ &$2k+1$ &$\frac{(2k+1)!}{(k!)^2}$&$6\lfloor \frac{k}{2}\rfloor$\\

$S_1 \times (S_2 \wr S_k)$&& $2k+1$ & $\frac{(2k+1)!}{2^k k!}$ & \\
\end{tabular}
\caption{Multiplicity-free intransitive groups, continued \label{table:intrans2}}
\end{center}
\end{table}
\renewcommand{\arraystretch}{1}

\newpage
\section{Association Schemes for Groups with Index less than 900}
The eigenvalues for the association schemes were calculated using
Akihide Hanaki's GAP package ``Elementary functions for association
schemes on GAP''\cite{ASfcn}.

$AGL(1,5) \cap A_5$ in $S_5$

\begin{displaymath}
\left(\begin{tabular}{cccc|c}
1 & 5 & 5 & 1 & 1 \\
1 & -5 & 5 & -1 & 1 \\
1 & -1 & -1 & 1 & 5 \\
1 & 1 & -1 & -1 & 5 \\
\end{tabular}
\right)
\end{displaymath}

$AGL(1,5)$ in $S_5$

\begin{displaymath}
\left(\begin{tabular}{cc|c}
1 & 5 & 1 \\
1 & -1 & 5 \\
\end{tabular}
\right)
\end{displaymath}

$PSL(2,5)$ in $S_6$

\begin{displaymath}
\left(\begin{tabular}{cccc|c}
1 & 5 & 5 & 1 & 1 \\
1 & -5 & 5 & -1 & 1 \\
1 & -1 & -1 & 1 & 5 \\
1 & 1 & -1 & -1 & 5 \\
\end{tabular}
\right)
\end{displaymath}

$PGL(2,5)$ in $S_6$

\begin{displaymath}
\left(\begin{tabular}{cc|c}
1 & 5 & 1 \\
1 & -1 & 5 \\
\end{tabular}
\right)
\end{displaymath}

$AGL(1,7)$ in $S_7$

\begin{displaymath}
\left(\begin{tabular}{ccccccc|c}
1 & 21 & 42 & 21 & 14 & 7 & 14 & 1 \\
1 & -6 & -3 & 3 & -1 & -2 & 8 & 14 \\
1 & 6 & -3 & -9 & -1 & 2 & 4 & 14 \\
1 & -7 & 2 & -3 & 6 & 3 & -2 & 15 \\
1 & 0 & 9 & -3 & -1 & -4 & -2 & 20 \\
1 & -1 & 2 & 3 & -6 & 3 & -2 & 21 \\
1 & 3 & -6 & 3 & 2 & -1 & -2 & 35 \\
\end{tabular}
\right)
\end{displaymath}

\newpage
$PSL(3,2)$ in $S_7$

\begin{displaymath}
\left(\begin{tabular}{cccc|c}
1 & 7 & 14 & 8 & 1 \\
1 & -7 & 14 & -8 & 1 \\
1 & -2 & -1 & 2 & 14 \\
1 & 2 & -1 & -2 & 14 \\
\end{tabular}
\right)
\end{displaymath}

$A\Gamma L(1,8)$ in $S_8$

\begin{displaymath}
\left(\begin{tabular}{cccccccc|c}
1 & 28 & 56 & 56 & 56 & 7 & 28 & 8 & 1 \\
1 & -28 & 56 & 56 & -56 & 7 & -28 & -8 & 1 \\
1 & -8 & -4 & -4 & 14 & 7 & -8 & 2 & 14 \\
1 & 8 & -4 & -4 & -14 & 7 & 8 & -2 & 14 \\
1 & -4 & 8 & -8 & -4 & -1 & 4 & 4 & 35 \\
1 & 4 & 8 & -8 & 4 & -1 & -4 & -4 & 35 \\
1 & -4 & -4 & 4 & 2 & -1 & 4 & -2 & 70 \\
1 & 4 & -4 & 4 & -2 & -1 & -4 & 2 & 70 \\
\end{tabular}
\right)
\end{displaymath}

$PGL(2,7)$ in $S_8$

\begin{displaymath}
\left(\begin{tabular}{ccccc|c}
1 & 28 & 56 & 14 & 21 & 1 \\
1 & -8 & -4 & 8 & 3 & 14 \\
1 & 8 & -4 & 4 & -9 & 14 \\
1 & -4 & 8 & -2 & -3 & 35 \\
1 & 2 & -4 & -2 & 3 & 56 \\
\end{tabular}
\right)
\end{displaymath}

$AGL(3,2)$ in $S_8$

\begin{displaymath}
\left(\begin{tabular}{cccc|c}
1 & 7 & 14 & 8 & 1 \\
1 & -7 & 14 & -8 & 1 \\
1 & -2 & -1 & 2 & 14 \\
1 & 2 & -1 & -2 & 14 \\
\end{tabular}
\right)
\end{displaymath}

\newpage
$AGL(2,3)$ in $S_9$

\begin{displaymath}
\left(\begin{tabular}{ccccccccc|c}
1 & 36 & 8 & 216 & 144 & 216 & 27 & 144 & 48 & 1 \\
1 & 0 & -4 & 36 & -12 & 0 & -9 & -24 & 12 & 42 \\
1 & 15 & 1 & 6 & 18 & -15 & 6 & -24 & -8 & 48 \\
1 & -9 & 5 & 0 & 18 & -27 & 0 & 0 & 12 & 56 \\
1 & -6 & 2 & 18 & -18 & 0 & 9 & 0 & -6 & 84 \\
1 & 6 & -2 & 6 & -6 & -24 & -3 & 24 & -2 & 84 \\
1 & 6 & 4 & -12 & -12 & 12 & -3 & 0 & 4 & 120 \\
1 & -4 & 0 & 0 & 8 & 8 & -5 & 0 & -8 & 189 \\
1 & -1 & -3 & -12 & 2 & 5 & 4 & 0 & 4 & 216 \\
\end{tabular}
\right)
\end{displaymath}

$P\Gamma L(2,8)$ in $S_9$

\begin{displaymath}
\left(\begin{tabular}{ccccc|c}
1 & 36 & 56 & 63 & 84 & 1 \\
1 & -36 & 56 & 63 & -84 & 1 \\
1 & 0 & 8 & -9 & 0 & 70 \\
1 & -6 & -4 & 3 & 6 & 84 \\
1 & 6 & -4 & 3 & -6 & 84 \\
\end{tabular}
\right)
\end{displaymath}

$S_2 \wr S_{2}$ in $S_{4}$

\begin{displaymath}
\left(\begin{tabular}{cc|c}
1 & 2 & 1 \\
1 & -1 & 2 \\
\end{tabular}
\right)
\end{displaymath}

$S_1 \times (S_2 \wr S_2)$ in $S_5$

\begin{displaymath}
\left(\begin{tabular}{cccc|c}
1 & 2 & 4 & 8 & 1 \\
1 & 2 & -1 & -2 & 4 \\
1 & -1 & -2 & 2 & 5 \\
1 & -1 & 2 & -2 & 5 \\
\end{tabular}
\right)
\end{displaymath}

$S_3 \wr S_{2}$ in $S_{6}$

\begin{displaymath}
\left(\begin{tabular}{cc|c}
1 & 9 & 1 \\
1 & -1 & 9 \\
\end{tabular}
\right)
\end{displaymath}

\newpage
$A_3 \wr S_2$ in $S_{6}$

\begin{displaymath}
\left(\begin{tabular}{cccccc|c}
1 & 2 & 1 & 9 & 18 & 9 & 1 \\
1 & -2 & 1 & -3 & 6 & -3 & 5 \\
1 & -2 & 1 & 3 & -6 & 3 & 5 \\
1 & 2 & 1 & -1 & -2 & -1 & 9 \\
1 & 0 & -1 & -3 & 0 & 3 & 10 \\
1 & 0 & -1 & 3 & 0 & -3 & 10 \\
\end{tabular}
\right)
\end{displaymath}

$(S_3 \wr S_2) \cap A_{6}$ in $S_{6}$

\begin{displaymath}
\left(\begin{tabular}{cccc|c}
1 & 1 & 9 & 9 & 1 \\
1 & -1 & -9 & 9 & 1 \\
1 & -1 & 1 & -1 & 9 \\
1 & 1 & -1 & -1 & 9 \\
\end{tabular}
\right)
\end{displaymath}

$S_1 \times (S_3 \wr S_2)$ in $S_7$

\begin{displaymath}
\left(\begin{tabular}{ccccc|c}
1 & 9 & 6 & 36 & 18 & 1 \\
1 & 9 & -1 & -6 & -3 & 6 \\
1 & -1 & 1 & 6 & -7 & 14 \\
1 & -1 & 4 & -6 & 2 & 14 \\
1 & -1 & -2 & 0 & 2 & 35 \\
\end{tabular}
\right)
\end{displaymath}

$S_1 \times ((S_3 \wr S_2)\cap A_{6})$ in $S_7$

\begin{displaymath}
\left(\begin{tabular}{cccccccccc|c}
1 & 1 & 9 & 9 & 6 & 6 & 36 & 36 & 18 & 18 & 1 \\
1 & -1 & -9 & 9 & -6 & 6 & 36 & -36 & -18 & 18 & 1 \\
1 & -1 & -9 & 9 & 1 & -1 & -6 & 6 & 3 & -3 & 6 \\
1 & 1 & 9 & 9 & -1 & -1 & -6 & -6 & -3 & -3 & 6 \\
1 & -1 & 1 & -1 & -4 & 4 & -6 & 6 & -2 & 2 & 14 \\
1 & -1 & 1 & -1 & -1 & 1 & 6 & -6 & 7 & -7 & 14 \\
1 & 1 & -1 & -1 & 1 & 1 & 6 & 6 & -7 & -7 & 14 \\
1 & 1 & -1 & -1 & 4 & 4 & -6 & -6 & 2 & 2 & 14 \\
1 & -1 & 1 & -1 & 2 & -2 & 0 & 0 & -2 & 2 & 35 \\
1 & 1 & -1 & -1 & -2 & -2 & 0 & 0 & 2 & 2 & 35 \\
\end{tabular}
\right)
\end{displaymath}

$S_4 \wr S_{2}$ in $S_{8}$

\begin{displaymath}
\left(\begin{tabular}{ccc|c}
1 & 16 & 18 & 1 \\
1 & -4 & 3 & 14 \\
1 & 2 & -3 & 20 \\
\end{tabular}
\right)
\end{displaymath}

$A_4 \wr S_2$ in $S_{8}$

\begin{displaymath}
\left(\begin{tabular}{cccccccc|c}
1 & 2 & 1 & 16 & 32 & 16 & 36 & 36 & 1 \\
1 & -2 & 1 & -16 & 32 & -16 & 36 & -36 & 1 \\
1 & -2 & 1 & 4 & -8 & 4 & 6 & -6 & 14 \\
1 & 2 & 1 & -4 & -8 & -4 & 6 & 6 & 14 \\
1 & -2 & 1 & -2 & 4 & -2 & -6 & 6 & 20 \\
1 & 2 & 1 & 2 & 4 & 2 & -6 & -6 & 20 \\
1 & 0 & -1 & -4 & 0 & 4 & 0 & 0 & 35 \\
1 & 0 & -1 & 4 & 0 & -4 & 0 & 0 & 35 \\
\end{tabular}
\right)
\end{displaymath}

$(S_4 \wr S_2) \cap A_{8}$ in $S_{8}$

\begin{displaymath}
\left(\begin{tabular}{cccccc|c}
1 & 1 & 16 & 16 & 18 & 18 & 1 \\
1 & -1 & -16 & 16 & 18 & -18 & 1 \\
1 & -1 & 4 & -4 & 3 & -3 & 14 \\
1 & 1 & -4 & -4 & 3 & 3 & 14 \\
1 & -1 & -2 & 2 & -3 & 3 & 20 \\
1 & 1 & 2 & 2 & -3 & -3 & 20 \\
\end{tabular}
\right)
\end{displaymath}

$S_1 \times (S_4 \wr S_2)$ in $S_9$

\begin{displaymath}
\left(\begin{tabular}{ccccccc|c}
1 & 16 & 18 & 8 & 96 & 144 & 32 & 1 \\
1 & 16 & 18 & -1 & -12 & -18 & -4 & 8 \\
1 & 2 & -3 & 6 & 2 & -18 & 10 & 27 \\
1 & -4 & 3 & 4 & -12 & 12 & -4 & 42 \\
1 & 2 & -3 & 1 & 12 & -3 & -10 & 48 \\
1 & -4 & 3 & -2 & 6 & -6 & 2 & 84 \\
1 & 2 & -3 & -2 & -6 & 6 & 2 & 105 \\
\end{tabular}
\right)
\end{displaymath}

\newpage
$S_1 \times ((S_4 \wr S_2)\cap A_{8})$ in $S_9$

\begin{displaymath}
\left(\begin{tabular}{cccccccccccccc|c}
1 & 1 & 16 & 16 & 18 & 18 & 8 & 8 & 96 & 96 & 144 & 144 & 32 & 32 & 1 \\
1 & -1 & -16 & 16 & 18 & -18 & -8 & 8 & 96 & -96 & -144 & 144 & 32 & -32 & 1 \\
1 & -1 & -16 & 16 & 18 & -18 & 1 & -1 & -12 & 12 & 18 & -18 & -4 & 4 & 8 \\
1 & 1 & 16 & 16 & 18 & 18 & -1 & -1 & -12 & -12 & -18 & -18 & -4 & -4 & 8 \\
1 & -1 & -2 & 2 & -3 & 3 & -6 & 6 & 2 & -2 & 18 & -18 & 10 & -10 & 27 \\
1 & 1 & 2 & 2 & -3 & -3 & 6 & 6 & 2 & 2 & -18 & -18 & 10 & 10 & 27 \\
1 & -1 & 4 & -4 & 3 & -3 & -4 & 4 & -12 & 12 & -12 & 12 & -4 & 4 & 42 \\
1 & 1 & -4 & -4 & 3 & 3 & 4 & 4 & -12 & -12 & 12 & 12 & -4 & -4 & 42 \\
1 & -1 & -2 & 2 & -3 & 3 & -1 & 1 & 12 & -12 & 3 & -3 & -10 & 10 & 48 \\
1 & 1 & 2 & 2 & -3 & -3 & 1 & 1 & 12 & 12 & -3 & -3 & -10 & -10 & 48 \\
1 & -1 & 4 & -4 & 3 & -3 & 2 & -2 & 6 & -6 & 6 & -6 & 2 & -2 & 84 \\
1 & 1 & -4 & -4 & 3 & 3 & -2 & -2 & 6 & 6 & -6 & -6 & 2 & 2 & 84 \\
1 & -1 & -2 & 2 & -3 & 3 & 2 & -2 & -6 & 6 & -6 & 6 & 2 & -2 & 105 \\
1 & 1 & 2 & 2 & -3 & -3 & -2 & -2 & -6 & -6 & 6 & 6 & 2 & 2 & 105 \\
\end{tabular}
\right)
\end{displaymath}

$S_5 \wr S_{2}$ in $S_{10}$

\begin{displaymath}
\left(\begin{tabular}{ccc|c}
1 & 25 & 100 & 1 \\
1 & 7 & -8 & 35 \\
1 & -3 & 2 & 90 \\
\end{tabular}
\right)
\end{displaymath}

$A_5 \wr S_2$ in $S_{10}$

\begin{displaymath}
\left(\begin{tabular}{cccccccc|c}
1 & 2 & 1 & 25 & 50 & 25 & 200 & 200 & 1 \\
1 & -2 & 1 & -15 & 30 & -15 & 40 & -40 & 9 \\
1 & 2 & 1 & 7 & 14 & 7 & -16 & -16 & 35 \\
1 & -2 & 1 & 5 & -10 & 5 & 20 & -20 & 42 \\
1 & -2 & 1 & -1 & 2 & -1 & -16 & 16 & 75 \\
1 & 2 & 1 & -3 & -6 & -3 & 4 & 4 & 90 \\
1 & 0 & -1 & -5 & 0 & 5 & 0 & 0 & 126 \\
1 & 0 & -1 & 5 & 0 & -5 & 0 & 0 & 126 \\
\end{tabular}
\right)
\end{displaymath}

\newpage
$(S_5 \wr S_2) \cap A_{10}$ in $S_{10}$

\begin{displaymath}
\left(\begin{tabular}{cccccc|c}
1 & 1 & 25 & 25 & 100 & 100 & 1 \\
1 & -1 & -25 & 25 & 100 & -100 & 1 \\
1 & -1 & -7 & 7 & -8 & 8 & 35 \\
1 & 1 & 7 & 7 & -8 & -8 & 35 \\
1 & -1 & 3 & -3 & 2 & -2 & 90 \\
1 & 1 & -3 & -3 & 2 & 2 & 90 \\
\end{tabular}
\right)
\end{displaymath}

$S_2 \wr S_{3}$ in $S_{6}$

\begin{displaymath}
\left(\begin{tabular}{ccc|c}
1 & 6 & 8 & 1 \\
1 & -3 & 2 & 5 \\
1 & 1 & -2 & 9 \\
\end{tabular}
\right)
\end{displaymath}

$S_2 \wr A_3$ in $S_{6}$

\begin{displaymath}
\left(\begin{tabular}{ccccc|c}
1 & 1 & 12 & 8 & 8 & 1 \\
1 & -1 & 0 & 4 & -4 & 5 \\
1 & 1 & -6 & 2 & 2 & 5 \\
1 & 1 & 2 & -2 & -2 & 9 \\
1 & -1 & 0 & -2 & 2 & 10 \\
\end{tabular}
\right)
\end{displaymath}

$(S_2 \wr S_3) \cap A_{6}$ in $S_{6}$

\begin{displaymath}
\left(\begin{tabular}{cccccc|c}
1 & 1 & 6 & 6 & 8 & 8 & 1 \\
1 & -1 & -6 & 6 & 8 & -8 & 1 \\
1 & -1 & 3 & -3 & 2 & -2 & 5 \\
1 & 1 & -3 & -3 & 2 & 2 & 5 \\
1 & -1 & -1 & 1 & -2 & 2 & 9 \\
1 & 1 & 1 & 1 & -2 & -2 & 9 \\
\end{tabular}
\right)
\end{displaymath}

$S_1 \times (S_2 \wr S_3)$ in $S_7$

\begin{displaymath}
\left(\begin{tabular}{ccccccc|c}
1 & 6 & 8 & 6 & 12 & 24 & 48 & 1 \\
1 & 6 & 8 & -1 & -2 & -4 & -8 & 6 \\
1 & -3 & 2 & -3 & 3 & 6 & -6 & 14 \\
1 & 1 & -2 & 1 & 7 & -6 & -2 & 14 \\
1 & 1 & -2 & 4 & -2 & 6 & -8 & 14 \\
1 & -3 & 2 & 2 & -2 & -4 & 4 & 21 \\
1 & 1 & -2 & -2 & -2 & 0 & 4 & 35 \\
\end{tabular}
\right)
\end{displaymath}

$S_3 \wr S_{3}$ in $S_{9}$

\begin{displaymath}
\left(\begin{tabular}{ccccc|c}
1 & 27 & 162 & 54 & 36 & 1 \\
1 & 11 & -6 & 6 & -12 & 27 \\
1 & 6 & -6 & -9 & 8 & 48 \\
1 & -3 & 12 & -6 & -4 & 84 \\
1 & -3 & -6 & 6 & 2 & 120 \\
\end{tabular}
\right)
\end{displaymath}

$(S_3 \wr S_3) \cap A_{9}$ in $S_{9}$

\begin{displaymath}
\left(\begin{tabular}{cccccccccc|c}
1 & 1 & 27 & 27 & 162 & 162 & 54 & 54 & 36 & 36 & 1 \\
1 & -1 & -27 & 27 & 162 & -162 & -54 & 54 & -36 & 36 & 1 \\
1 & -1 & -11 & 11 & -6 & 6 & -6 & 6 & 12 & -12 & 27 \\
1 & 1 & 11 & 11 & -6 & -6 & 6 & 6 & -12 & -12 & 27 \\
1 & -1 & -6 & 6 & -6 & 6 & 9 & -9 & -8 & 8 & 48 \\
1 & 1 & 6 & 6 & -6 & -6 & -9 & -9 & 8 & 8 & 48 \\
1 & -1 & 3 & -3 & 12 & -12 & 6 & -6 & 4 & -4 & 84 \\
1 & 1 & -3 & -3 & 12 & 12 & -6 & -6 & -4 & -4 & 84 \\
1 & -1 & 3 & -3 & -6 & 6 & -6 & 6 & -2 & 2 & 120 \\
1 & 1 & -3 & -3 & -6 & -6 & 6 & 6 & 2 & 2 & 120 \\
\end{tabular}
\right)
\end{displaymath}

$S_2 \wr S_{4}$ in $S_{8}$

\begin{displaymath}
\left(\begin{tabular}{ccccc|c}
1 & 12 & 32 & 12 & 48 & 1 \\
1 & -6 & 8 & 3 & -6 & 14 \\
1 & 2 & -8 & 7 & -2 & 14 \\
1 & 5 & 4 & -2 & -8 & 20 \\
1 & -1 & -2 & -2 & 4 & 56 \\
\end{tabular}
\right)
\end{displaymath}

$S_2 \wr A_4$ in $S_{8}$

\begin{displaymath}
\left(\begin{tabular}{ccccccc|c}
1 & 1 & 24 & 32 & 32 & 24 & 96 & 1 \\
1 & 1 & -12 & 8 & 8 & 6 & -12 & 14 \\
1 & 1 & 4 & -8 & -8 & 14 & -4 & 14 \\
1 & 1 & 10 & 4 & 4 & -4 & -16 & 20 \\
1 & -1 & 0 & -8 & 8 & 0 & 0 & 35 \\
1 & 1 & -2 & -2 & -2 & -4 & 8 & 56 \\
1 & -1 & 0 & 4 & -4 & 0 & 0 & 70 \\
\end{tabular}
\right)
\end{displaymath}

\newpage
$S_1 \times (S_2 \wr S_4)$ in $S_9$

\begin{displaymath}
\left(\begin{tabular}{cccccccccccc|c}
1 & 12 & 32 & 12 & 48 & 8 & 48 & 64 & 48 & 96 & 192 & 384 & 1 \\
1 & 12 & 32 & 12 & 48 & -1 & -6 & -8 & -6 & -12 & -24 & -48 & 8 \\
1 & 5 & 4 & -2 & -8 & 6 & 8 & -8 & 22 & -12 & 32 & -48 & 27 \\
1 & -6 & 8 & 3 & -6 & -4 & 12 & -8 & 12 & -12 & -24 & 24 & 42 \\
1 & 2 & -8 & 7 & -2 & 4 & 4 & -8 & 4 & 28 & -24 & -8 & 42 \\
1 & 5 & 4 & -2 & -8 & 1 & 13 & 22 & -8 & -2 & -18 & -8 & 48 \\
1 & -6 & 8 & 3 & -6 & 2 & -6 & 4 & -6 & 6 & 12 & -12 & 84 \\
1 & 2 & -8 & 7 & -2 & -2 & -2 & 4 & -2 & -14 & 12 & 4 & 84 \\
1 & 5 & 4 & -2 & -8 & -2 & -8 & -8 & -2 & 4 & 0 & 16 & 105 \\
1 & -1 & -2 & -2 & 4 & 4 & -8 & 4 & 4 & -8 & -12 & 16 & 120 \\
1 & -1 & -2 & -2 & 4 & 1 & 7 & -8 & -8 & -2 & 6 & 4 & 168 \\
1 & -1 & -2 & -2 & 4 & -3 & -1 & 4 & 4 & 6 & 2 & -12 & 216 \\
\end{tabular}
\right)
\end{displaymath}

$(S_3^2 \cap A_{6}) \rtimes S_2$ in $S_{6}$

\begin{displaymath}
\left(\begin{tabular}{cccc|c}
1 & 1 & 9 & 9 & 1 \\
1 & -1 & -3 & 3 & 5 \\
1 & -1 & 3 & -3 & 5 \\
1 & 1 & -1 & -1 & 9 \\
\end{tabular}
\right)
\end{displaymath}

$(S_4^2 \cap A_{8} )\rtimes S_2$ in $S_{8}$

\begin{displaymath}
\left(\begin{tabular}{cccccc|c}
1 & 1 & 16 & 16 & 18 & 18 & 1 \\
1 & -1 & -16 & 16 & 18 & -18 & 1 \\
1 & -1 & 4 & -4 & 3 & -3 & 14 \\
1 & 1 & -4 & -4 & 3 & 3 & 14 \\
1 & -1 & -2 & 2 & -3 & 3 & 20 \\
1 & 1 & 2 & 2 & -3 & -3 & 20 \\
\end{tabular}
\right)
\end{displaymath}

$(S_5^2 \cap A_{10}) \rtimes S_2$ in $S_{10}$

\begin{displaymath}
\left(\begin{tabular}{cccccc|c}
1 & 1 & 25 & 25 & 100 & 100 & 1 \\
1 & -1 & -15 & 15 & 20 & -20 & 9 \\
1 & 1 & 7 & 7 & -8 & -8 & 35 \\
1 & -1 & 5 & -5 & 10 & -10 & 42 \\
1 & -1 & -1 & 1 & -8 & 8 & 75 \\
1 & 1 & -3 & -3 & 2 & 2 & 90 \\
\end{tabular}
\right)
\end{displaymath}

$(S_3^3 \cap A_{9}) \rtimes S_3$ in $S_{9}$

\begin{displaymath}
\left(\begin{tabular}{ccccccccc|c}
1 & 1 & 27 & 27 & 162 & 162 & 54 & 54 & 72 & 1 \\
1 & 1 & 11 & 11 & -6 & -6 & 6 & 6 & -24 & 27 \\
1 & -1 & -9 & 9 & 0 & 0 & 18 & -18 & 0 & 28 \\
1 & -1 & 9 & -9 & 18 & -18 & 6 & -6 & 0 & 42 \\
1 & -1 & -6 & 6 & 18 & -18 & -9 & 9 & 0 & 48 \\
1 & 1 & 6 & 6 & -6 & -6 & -9 & -9 & 16 & 48 \\
1 & 1 & -3 & -3 & 12 & 12 & -6 & -6 & -8 & 84 \\
1 & 1 & -3 & -3 & -6 & -6 & 6 & 6 & 4 & 120 \\
1 & -1 & 1 & -1 & -10 & 10 & -2 & 2 & 0 & 162 \\
\end{tabular}
\right)
\end{displaymath}

$S_1 \times ( (S_3^2 \cap A_{6}) \rtimes S_2)$ in $S_{7}$

\begin{displaymath}
\left(\begin{tabular}{cccccccccc|c}
1 & 1 & 9 & 9 & 6 & 6 & 36 & 36 & 18 & 18 & 1 \\
1 & -1 & -3 & 3 & -5 & 5 & 6 & -6 & 9 & -9 & 6 \\
1 & 1 & 9 & 9 & -1 & -1 & -6 & -6 & -3 & -3 & 6 \\
1 & -1 & -3 & 3 & 0 & 0 & 6 & -6 & -6 & 6 & 14 \\
1 & -1 & 3 & -3 & -3 & 3 & -6 & 6 & -3 & 3 & 14 \\
1 & 1 & -1 & -1 & 1 & 1 & 6 & 6 & -7 & -7 & 14 \\
1 & 1 & -1 & -1 & 4 & 4 & -6 & -6 & 2 & 2 & 14 \\
1 & -1 & -3 & 3 & 2 & -2 & -8 & 8 & 2 & -2 & 15 \\
1 & -1 & 3 & -3 & 2 & -2 & 4 & -4 & 2 & -2 & 21 \\
1 & 1 & -1 & -1 & -2 & -2 & 0 & 0 & 2 & 2 & 35 \\
\end{tabular}
\right)
\end{displaymath}

\newpage
$S_1 \times ( (S_4^2 \cap A_{8} )\rtimes S_2)$ in $S_{9}$

\begin{displaymath}
\left(\begin{tabular}{cccccccccccccc|c}
1 & 1 & 16 & 16 & 18 & 18 & 8 & 8 & 96 & 96 & 144 & 144 & 32 & 32 & 1 \\
1 & -1 & -16 & 16 & 18 & -18 & -8 & 8 & 96 & -96 & -144 & 144 & 32 & -32 & 1 \\
1 & -1 & -16 & 16 & 18 & -18 & 1 & -1 & -12 & 12 & 18 & -18 & -4 & 4 & 8 \\
1 & 1 & 16 & 16 & 18 & 18 & -1 & -1 & -12 & -12 & -18 & -18 & -4 & -4 & 8 \\
1 & -1 & -2 & 2 & -3 & 3 & -6 & 6 & 2 & -2 & 18 & -18 & 10 & -10 & 27 \\
1 & 1 & 2 & 2 & -3 & -3 & 6 & 6 & 2 & 2 & -18 & -18 & 10 & 10 & 27 \\
1 & -1 & 4 & -4 & 3 & -3 & -4 & 4 & -12 & 12 & -12 & 12 & -4 & 4 & 42 \\
1 & 1 & -4 & -4 & 3 & 3 & 4 & 4 & -12 & -12 & 12 & 12 & -4 & -4 & 42 \\
1 & -1 & -2 & 2 & -3 & 3 & -1 & 1 & 12 & -12 & 3 & -3 & -10 & 10 & 48 \\
1 & 1 & 2 & 2 & -3 & -3 & 1 & 1 & 12 & 12 & -3 & -3 & -10 & -10 & 48 \\
1 & -1 & 4 & -4 & 3 & -3 & 2 & -2 & 6 & -6 & 6 & -6 & 2 & -2 & 84 \\
1 & 1 & -4 & -4 & 3 & 3 & -2 & -2 & 6 & 6 & -6 & -6 & 2 & 2 & 84 \\
1 & -1 & -2 & 2 & -3 & 3 & 2 & -2 & -6 & 6 & -6 & 6 & 2 & -2 & 105 \\
1 & 1 & 2 & 2 & -3 & -3 & -2 & -2 & -6 & -6 & 6 & 6 & 2 & 2 & 105 \\
\end{tabular}
\right)
\end{displaymath}

$S_1 \times AGL(1,5)$ in $S_6$
\begin{displaymath}
\left(\begin{tabular}{ccccc|c}
1 & 5 & 5 & 20 & 5 & 1 \\
1 & -1 & -1 & -4 & 5 & 5 \\
1 & 5 & -1 & -4 & -1 & 5 \\
1 & -1 & -3 & 4 & -1 & 9 \\
1 & -1 & 2 & -1 & -1 & 16 \\
\end{tabular}
\right)
\end{displaymath}

$S_1 \times PSL(2,5)$ in $S_7$
\begin{displaymath}
\left(\begin{tabular}{cccccccc|c}
1 & 5 & 5 & 1 & 6 & 30 & 30 & 6 & 1 \\
1 & -5 & 5 & -1 & -6 & 30 & -30 & 6 & 1 \\
1 & -5 & 5 & -1 & 1 & -5 & 5 & -1 & 6 \\
1 & 5 & 5 & 1 & -1 & -5 & -5 & -1 & 6 \\
1 & -1 & -1 & 1 & -3 & 3 & 3 & -3 & 14 \\
1 & 1 & -1 & -1 & 3 & 3 & -3 & -3 & 14 \\
1 & -1 & -1 & 1 & 2 & -2 & -2 & 2 & 21 \\
1 & 1 & -1 & -1 & -2 & -2 & 2 & 2 & 21 \\
\end{tabular}
\right)
\end{displaymath}

$S_1 \times PGL(2,5)$ in $S_7$
\begin{displaymath}
\left(\begin{tabular}{cccc|c}
1 & 5 & 6 & 30 & 1 \\
1 & 5 & -1 & -5 & 6 \\
1 & -1 & -3 & 3 & 14 \\
1 & -1 & 2 & -2 & 21 \\
\end{tabular}
\right)
\end{displaymath}

$S_1 \times PSL(3,2)$ in $S_8$
\begin{displaymath}
\left(\begin{tabular}{cccccccccc|c}
1 & 7 & 14 & 8 & 7 & 42 & 7 & 42 & 56 & 56 & 1 \\
1 & -7 & 14 & -8 & -7 & 42 & 7 & -42 & -56 & 56 & 1 \\
1 & -7 & 14 & -8 & 1 & -6 & -1 & 6 & 8 & -8 & 7 \\
1 & 7 & 14 & 8 & -1 & -6 & -1 & -6 & -8 & -8 & 7 \\
1 & -2 & -1 & 2 & -2 & -3 & 7 & -12 & 14 & -4 & 14 \\
1 & 2 & -1 & -2 & 2 & -3 & 7 & 12 & -14 & -4 & 14 \\
1 & -2 & -1 & 2 & -4 & 9 & -1 & 6 & -2 & -8 & 28 \\
1 & 2 & -1 & -2 & 4 & 9 & -1 & -6 & 2 & -8 & 28 \\
1 & -2 & -1 & 2 & 2 & -3 & -1 & 0 & -2 & 4 & 70 \\
1 & 2 & -1 & -2 & -2 & -3 & -1 & 0 & 2 & 4 & 70 \\
\end{tabular}
\right)
\end{displaymath}

$S_1 \times AGL(3,2)$ in $S_9$
\begin{displaymath}
\left(\begin{tabular}{cccccccc|c}
1 & 7 & 14 & 8 & 8 & 56 & 112 & 64 & 1 \\
1 & -7 & 14 & -8 & -8 & 56 & -112 & 64 & 1 \\
1 & -7 & 14 & -8 & 1 & -7 & 14 & -8 & 8 \\
1 & 7 & 14 & 8 & -1 & -7 & -14 & -8 & 8 \\
1 & -2 & -1 & 2 & -4 & 8 & 4 & -8 & 42 \\
1 & 2 & -1 & -2 & 4 & 8 & -4 & -8 & 42 \\
1 & -2 & -1 & 2 & 2 & -4 & -2 & 4 & 84 \\
1 & 2 & -1 & -2 & -2 & -4 & 2 & 4 & 84 \\
\end{tabular}
\right)
\end{displaymath}

$S_2 \times AGL(1,5)$ in $S_7$
\begin{displaymath}
\left(\begin{tabular}{ccccccc|c}
1 & 5 & 10 & 40 & 10 & 20 & 40 & 1 \\
1 & 5 & 3 & 12 & 3 & -8 & -16 & 6 \\
1 & -1 & -5 & 4 & 1 & 8 & -8 & 14 \\
1 & 5 & -2 & -8 & -2 & 2 & 4 & 14 \\
1 & -1 & 0 & -6 & 6 & -2 & 2 & 21 \\
1 & -1 & -2 & 4 & -2 & -4 & 4 & 35 \\
1 & -1 & 4 & -2 & -2 & 2 & -2 & 35 \\
\end{tabular}
\right)
\end{displaymath}

$S_2 \times PSL(2,5)$ in $S_8$
\begin{displaymath}
\left(\begin{tabular}{ccccccccccc|c}
1 & 5 & 5 & 1 & 12 & 60 & 60 & 12 & 30 & 120 & 30 & 1 \\
1 & -5 & 5 & -1 & -6 & 30 & -30 & 6 & 0 & 0 & 0 & 7 \\
1 & 5 & 5 & 1 & 4 & 20 & 20 & 4 & -10 & -40 & -10 & 7 \\
1 & -1 & -1 & 1 & -6 & 6 & 6 & -6 & 12 & -24 & 12 & 14 \\
1 & 5 & 5 & 1 & -2 & -10 & -10 & -2 & 2 & 8 & 2 & 20 \\
1 & -5 & 5 & -1 & 2 & -10 & 10 & -2 & 0 & 0 & 0 & 21 \\
1 & 1 & -1 & -1 & 6 & 6 & -6 & -6 & 6 & 0 & -6 & 28 \\
1 & 1 & -1 & -1 & -4 & -4 & 4 & 4 & 6 & 0 & -6 & 42 \\
1 & -1 & -1 & 1 & 4 & -4 & -4 & 4 & 2 & -4 & 2 & 56 \\
1 & -1 & -1 & 1 & -2 & 2 & 2 & -2 & -4 & 8 & -4 & 70 \\
1 & 1 & -1 & -1 & 0 & 0 & 0 & 0 & -6 & 0 & 6 & 70 \\
\end{tabular}
\right)
\end{displaymath}

$S_2 \times PGL(2,5)$ in $S_8$
\begin{displaymath}
\left(\begin{tabular}{cccccc|c}
1 & 5 & 12 & 60 & 30 & 60 & 1 \\
1 & 5 & 4 & 20 & -10 & -20 & 7 \\
1 & -1 & -6 & 6 & 12 & -12 & 14 \\
1 & 5 & -2 & -10 & 2 & 4 & 20 \\
1 & -1 & 4 & -4 & 2 & -2 & 56 \\
1 & -1 & -2 & 2 & -4 & 4 & 70 \\
\end{tabular}
\right)
\end{displaymath}

$S_3 \times AGL(1,5)$ in $S_8$
\begin{displaymath}
\left(\begin{tabular}{cccccccc|c}
1 & 5 & 15 & 60 & 15 & 60 & 120 & 60 & 1 \\
1 & 5 & 7 & 28 & 7 & -4 & -8 & -36 & 7 \\
1 & 5 & 1 & 4 & 1 & -10 & -20 & 18 & 20 \\
1 & 5 & -3 & -12 & -3 & 6 & 12 & -6 & 28 \\
1 & -1 & 1 & -8 & 7 & -4 & 4 & 0 & 56 \\
1 & -1 & 6 & -3 & -3 & 6 & -6 & 0 & 64 \\
1 & -1 & -5 & 4 & 1 & 8 & -8 & 0 & 70 \\
1 & -1 & -1 & 4 & -3 & -8 & 8 & 0 & 90 \\
\end{tabular}
\right)
\end{displaymath}

$S_3 \times PGL(2,5)$ in $S_9$
\begin{displaymath}
\left(\begin{tabular}{ccccccc|c}
1 & 5 & 18 & 90 & 90 & 180 & 120 & 1 \\
1 & 5 & 9 & 45 & 0 & 0 & -60 & 8 \\
1 & 5 & 2 & 10 & -14 & -28 & 24 & 27 \\
1 & 5 & -3 & -15 & 6 & 12 & -6 & 48 \\
1 & -1 & -6 & 6 & 12 & -12 & 0 & 84 \\
1 & -1 & 6 & -6 & 6 & -6 & 0 & 120 \\
1 & -1 & -1 & 1 & -8 & 8 & 0 & 216 \\
\end{tabular}
\right)
\end{displaymath}

$S_4 \times AGL(1,5)$ in $S_9$
\begin{displaymath}
\left(\begin{tabular}{ccccccccc|c}
1 & 5 & 20 & 80 & 20 & 120 & 240 & 240 & 30 & 1 \\
1 & 5 & 11 & 44 & 11 & 12 & 24 & -84 & -24 & 8 \\
1 & 5 & 4 & 16 & 4 & -16 & -32 & 0 & 18 & 27 \\
1 & 5 & -4 & -16 & -4 & 12 & 24 & -24 & 6 & 42 \\
1 & 5 & -1 & -4 & -1 & -6 & -12 & 30 & -12 & 48 \\
1 & -1 & 8 & -4 & -4 & 12 & -12 & 0 & 0 & 105 \\
1 & -1 & 2 & -10 & 8 & -6 & 6 & 0 & 0 & 120 \\
1 & -1 & 0 & 4 & -4 & -12 & 12 & 0 & 0 & 189 \\
1 & -1 & -5 & 4 & 1 & 8 & -8 & 0 & 0 & 216 \\
\end{tabular}
\right)
\end{displaymath}

\end{document}